\documentclass[12pt]{amsart}
\usepackage{latexsym, amsmath,amssymb}
\usepackage{hyperref}
\usepackage{xcolor}
\hypersetup{
    colorlinks,
    linkcolor={red!50!black},
    citecolor={blue!50!black},
    urlcolor={blue!80!black}
}

\usepackage{graphicx}
\usepackage{float}

\providecommand{\abs}[1]{\lvert#1\rvert}
\providecommand{\norm}[1]{\lVert#1\rVert} 
 \numberwithin{equation}{section}

\newenvironment{claim}{\medskip\noindent{\it Claim:\/} }{\medskip}

\def\XXint#1#2#3{{\setbox0=\hbox{$#1{#2#3}{%
\int}$ }
\vcenter{\hbox{$#2#3$ }}\kern-.6\wd0}}

\renewcommand{\epsilon}{\varepsilon}
\newtheorem{theorem}{Theorem}

\newtheorem{lemma}[theorem]{Lemma}
\newtheorem{corr}[theorem]{Corollary}

\newtheorem{proposition}[theorem]{Proposition}
\newtheorem{deff}[theorem]{Definition}

\newcommand{\bth}{\begin{theorem}}
\newcommand{\ble}{\begin{lemma}}
\newcommand{\bcor}{\begin{corr}}

\newcommand{\bdeff}{\begin{deff}}

\newcommand{\bprop}{\begin{proposition}}
\newcommand{\ele}{\end{lemma}}
\newcommand{\ecor}{\end{corr}}
\newcommand{\edeff}{\end{deff}}

\numberwithin{theorem}{section}

\newcommand{\eprop}{\end{proposition}}

\renewcommand{\Pi}{\varPi}

\renewcommand{\epsilon}{\varepsilon}

\begin{document}

\title[The critical quasigeostrophic equation on $\mathbb{S}^2$]
{Continuity of weak solutions of\\ the critical quasigeostrophic equation on $\mathbb{S}^2$}
\author[D. Alonso-Or\'an]{Diego Alonso-Or\'an}
\address{Instituto de Ciencias Matem\'aticas CSIC-UAM-UC3M-UCM -- Departamento de Matem\'aticas (Universidad Aut\'onoma de Madrid), 28049 Madrid, Spain} 
\email{diego.alonso@icmat.es}
\author[A. C\'ordoba]{Antonio C\'ordoba}
\address{Instituto de Ciencias Matem\'aticas CSIC-UAM-UC3M-UCM -- Departamento de Matem\'aticas (Universidad Aut\'onoma de Madrid), 28049 Madrid, Spain} 
\email{antonio.cordoba@uam.es}
\author[A. D. Mart\'inez]{\'Angel D. Mart\'inez}
\address{Instituto de Ciencias Matem\'aticas (CSIC-UAM-UC3M-UCM) -- Departamento de Matem\'aticas (Universidad Aut\'onoma de Madrid), 28049 Madrid, Spain} 
\email{angel.martinez@icmat.es}

\begin{abstract}
In this paper we provide regularity results for active scalars that are weak solutions of almost critical drift-diffusion equations in general surfaces. This includes models of anisotropic non-homogeneous media and the physically motivated case of the two-dimensional sphere. Our finest result deals with the critical surface quasigeostrophic equation on the round sphere.
\end{abstract}

\maketitle

\section{\textbf{Introduction}}	

General drift-diffusion equations refer to evolution equations of the form
\[\partial_t\theta+u\cdot\nabla\theta+\kappa\Lambda^{\alpha}\theta=0\]
where $u$ is a vector field, $\alpha\in(0,2)$ is a difussion exponent and $\kappa$ a positive constant. This type of equations have been studied intensively by a number of authors during the last decades. We are mainly interested in the case where the drift velocity is given by $u=\nabla^{\perp}\Lambda^{-1}\theta$ which corresponds to the surface quasigeostrophic equation, which describes the evolution of a temperature field in a rapidly rotating stratified fluid with potential vorticity \cite{held}. The numerical and analytical study of the surface quasigeostrophic equation started in \cite{CMT}, motivated also its analogy with the three-dimensional Euler equation given in vorticity form. A remarkable singularity scenario was dismissed by D. C\'ordoba in \cite{Cdie}, but whether or not solutions to the inviscid surface quasigeostrophic equation, $\kappa=0$, develop singularities in finite time, represents a major open problem.

The issue of global regularity versus finite time blow up for the fractional surface quasigeostrophic equation  has attracted a lot of attention. The subcritical case ($\alpha>1$) is well understood and essentialy solved in \cite{R,CW}. Global regularity in the critical case ($\alpha=1$) is quite more challenging due to the possible balance between opposite strengths of the nonlinear and the dissipative term. Therefore perturbative methods are not useful anymore and more refined techniques come in to play.  Constantin, C\'ordoba and Wu adressed for the first time the global regularity for the critical case under a small data hypothesis and global existence was obtained (cf. \cite{CCW}). Two different approaches developed independently, one by Kiselev, Nazarov and Volberg and the other by Caffareli and Vasseur, to tackle this sophisticated problem (cf. \cite{KNV, CV}). Then a refinement of the pointwise inequality of C\'ordoba-C\'ordoba (cf. \cite{CC, CC2}) also proved its usefulness for such matter in the work of Constantin and Vicol, who gave a third proof of global regularity which relies on non linear lower bounds for the fractional Laplacian (cf. \cite{CC, CVi}). Recent works by Constantin and Ignatova  extend this results for bounded domains (cf. \cite{CI1,CI2}). To the best of our knowledge global regularity for the supercritical case remains unsolved. 

In this paper we address the same equation considered in a compact orientable surface $M$ with riemannian metric $g$, where $\Lambda$ will denote the square root of its Laplace-Beltrami operator $-\Delta_g$. In the particular case of the two-dimensional round sphere we get the following:

\begin{theorem}\label{I}
Let $\theta_0\in L^2(\mathbb{S}^2)$ and $\theta$ a weak solution of the following Cauchy problem
\[\left\{\begin{array}{l}
\partial_t\theta+u\cdot\nabla_g\theta=-\Lambda\theta\\
\theta(x,0)=\theta_0
\end{array}\right.\]
where $u=\nabla_g^{\perp}\Lambda^{-1}\theta$. Then $\theta(x,t)$ is continuous for every $t>0$.
\end{theorem}

In the particular case of the two dimensional sphere an explicit computation shows $\textrm{div}_g\nabla_g^{\perp}=0$. However this fact also holds on any two dimensional Riemannian manifold. The higher dimensional analogue is more delicate, but in even dimensions and in the presence of a symplectic structure, there is a canonical construction of an orthogonal gradient of a function $f$ such that its divergence vanishes; that is, the vector field from the statement is incompressible. Therefore it also makes sense to study the surface quasigeostrophic equation there.

Interestingly enough our proof of the stated theorem breaks down for higher dimensional spheres. However it shows a stronger result, namely, it provides an explicit modulus of continuity, following the non local version of De Giorgi's robust strategy introduced by Caffarelli and Vasseur (see \cite{CV, DG}; cf. \cite{N}). Their results is claimed for drift-diffusion equations whose divergence free velocity field $u$ has bounded mean oscillation. However, our situation is not that fortunate because curvature matters and some extra hypothesis should be made in order to extend their analysis. We isolate part of that fact as a separate theorem which holds for compact riemannian manifolds.

\begin{theorem}\label{II}
Let $\theta_0\in L^2(M)$ and $\theta$ a weak solution of the following Cauchy problem
\begin{eqnarray}\label{sqg}
\left\{\begin{array}{l}
\partial_t\theta+u\cdot\nabla_g\theta=-\Lambda\theta\\ 
\theta(x,0)=\theta_0 \nonumber
\end{array}\right.
\end{eqnarray}
where the divergence free velocity field $u\in L^{\infty}(M)$ uniformly in time. Then $\theta(x,t)$ is of class $C^{\alpha}$ for any $t>0$.
\end{theorem}

The hypothesis is expected to be satisfied in the subcritical regime, where standard harmonic analysis techniques prove it when $u=\nabla_g^{\perp}\Lambda^{-1-\epsilon}\theta$ for $\epsilon>0$, but we shall not consider here that case. In the more interesting critical case $\epsilon=0$ certain difficulties arise coming from the non locality, which are responsible for the limitations appearing in the statement of Theorem \ref{I}. The main one to be overcomed is the pointwise inequality in this setting (cf. \cite{CM}) and the existence of adequate barriers at different scales. In section \S\ref{barrier} we prove a quantitative maximum principle for certain family of barriers adapted to the local geometry of the manifold. Since there is no device for rescaling the equation, as opposed to the euclidean setting, our statements have to made special emphasis to take into account the scale influence in the arguments. The rest of the paper presents the proof in full details.

\section{\textbf{Scheme of the proof (d'apr\'es Caffarelli-Vasseur)}}\label{cafvas}

For the sake of completeness, in the appendix ~\ref{debiles} we provide a proof of the existence of global weak solutions to the Cauchy problem for the critical quasigeostrophic equation. We need to use a fractional Sobolev embedding on compact manifolds which is not easy to find in the literature, but a detailed proof of that fact was included in a recent paper of the authors \cite{AOMC} (cf. \cite{Au}). 

The paper follows to some extent the structure of the work of Caffarelli and Vasseur \cite{CV}. In section \S\ref{infbound} we prove an uniform bound for the essential supremum of a global weak solution in space and strictly positive time $t\geq t_0>0$. The bound depends on $t_0$ and the initial energy. The proof is based upon a non linear inequality as in De Giorgi's iteration scheme. Once this is achieved, we may treat the problem as if it was linear, forgetting the $\theta$ dependence of $u$. Therefore, we get a drift diffusion equation where the drift is in some appropiate functional space, which will be enough to prove the local energy estimate, instrumental to control the oscillation decay. Let us digress briefly why a control on the oscillation decay is a convenient strategy: the oscillation of a function $f$ in a ball $B$ is defined as
\[\textrm{osc}_{B}f=\sup_{x,y\in B}|f(x)-f(y)|.\]
If one has an estimate of the type
\[\textrm{osc}_{B_{h/2}}f\leq\delta\cdot\textrm{osc}_{B_h}f\]
for some fixed $\delta<1$ and any scale $h$, then it is easy to prove $f\in C^{\alpha}$ for some $\alpha$ which depends on $\delta$. Indeed, we just need to control the norm
\[\|f\|_{C^{\alpha}}=\|f\|_{L^{\infty}}+\sup_{x,y}\frac{|f(x)-f(y)|}{d(x,y)^{\alpha}}.\]
The $L^{\infty}$ norm will be controlled in the first part of the proof, see \S\ref{infbound}. For the second term let us observe that
\[\frac{|f(x)-f(y)|}{d(x,y)^{\alpha}}\leq\sup_k\sup_{y\in A_k(x)}|f(x)-f(y)|2^{\alpha k}\]
where $A_k(x)=B_{1/2^{k-1}}(x)\setminus B_{1/2^k}(x)$ are coronas centered in $x$. But the right hand side is bounded by
\[\sup_{k}2^{\alpha k}\textrm{osc}_{B_{1/2^{k-1}}}f \leq\sup_{k} 2^{\alpha k} \delta^k\textrm{osc}_{B_{1}}f\leq 2\|f\|_{L^{\infty}}(2^{\alpha}\delta)^k.\]
Now it becomes obvious that choosing $\alpha<\log_2(\delta)$ would be enough to obtain the bound. Notice that we used the notation $d(x,y)$ instead of $|x-y|$ since we are not dealing with the euclidean metric. A version of the oscillation decay will be achieved using De Giorgi's iteration scheme if the initial energy is small in \S\ref{small}. It is at this stage where the precise control on barrier functions developed in \S\ref{barrier} becomes crucial. This together with the local energy inequality proved in \S\ref{localenergy} yields a non linear energy inequality which is at the heart of the oscillation decay. Finally, a non-local version of De Giorgi's isoperimetric inequality will drop the small mean energy condition. This is done in \S\ref{big} which, though quite similar, differs from the original treatment \cite{CV} where the argument has been carefully adapted to respect different scales. This would prove Theorem \ref{II}. 

Theorem \ref{I} follows from the latter but not so inmediately since interpolation arguments might lead to a $L^{2n}(M)$ estimate which does not imply the needed scale decay in small balls. In order to handle this problem we introduce an auxiliary function satisfying the same equation and whose $L^{2n}(B_g(h))$ norm is controled. The arguments leading to Theorem \ref{II} will then be applicable but we cannot obtain the necessary decay to imply H\"older continuity. However it will be enough to establish a modulus of continuity of the form $\omega(\rho)=\log(1/\rho)^{-\alpha}$ for some $\alpha>0$. It should be remarked that, unfortunately, it falls close, but not enough, to satisfy the classical Dini condition under which Theorem \ref{II} would be applicable.

\section{\textbf{$L^{\infty}_{x,t}$ bound}}\label{infbound}

In this section we illustrate De Giorgi's method which will be based on a non linear inequality for some sort of energy. But a more subtle needed version will be exposed in sections \ref{small} and \ref{big}.

\begin{proposition}\label{l2linf}
Let $\theta(x,t)$ be a weak solution (cf. Appendix \ref{debiles}). Then for any fixed $t_{0}>0$ there exists a positive constant $C$ that will depend on $\theta_0$'s energy, $t_0$ and the manifold such that 
\[\abs{\theta(x,t)}\leq C\textrm{ for any $x\in M$ and any $t>t_{0}$.}\]
\end{proposition}

\textit{Remark:}\label{constantes} in the rest of the paper all constants $C$ will be assumed to depend implicitly on quantities that are considered to be constant. In particular, they will have to be scale independent. Notice also that the constant might differ from line to line for the sake of the exposition's clearness.

\textsc{Proof of Proposition \ref{l2linf}	:} we will proceed using a nonlinear energy inequality for consecutive energy truncations which is based on the interplay between a global energy inequality and  Sobolev inequality. Let us assume without loss of generality that $\int_M\theta(x,t)d\textrm{vol}_g(x)=0$ and define the truncation levels as follows:
\[\ell_{k}=C(1-2^{-k})\]
where $C$ will be chosen later to be large enough. The $k$th truncation of $\theta$ at the level $\ell_k$ will be denoted by $\theta_{k}=(\theta-\ell_{k})_{+}$. Notice $(a)_+=\max\{a,0\}$ is a convex function. One can derive a differential inequality for the truncations using the C\'ordoba-C\'ordoba pointwise inequality for fractional powers of the Laplace-Beltrami operator on manifolds (cf. \cite{CM}) 
\[\partial_t\theta_k+u\cdot\nabla_g\theta_k\leq-\Lambda\theta_k\]
multiplying this by $\theta_k$, integrating in $M$, using that $u$ is divergence free and playing around with the fractional Laplace-Beltrami operator the following holds
\[\partial_{t}\int_{M}\theta^{2}_{k} d\textrm{vol}_{g}(x)+ \int_{M}\abs{\Lambda^{1/2}\theta_{k}}^{2} \ d\textrm{vol}_{g}(x) \leq 0.\]

Let us introduce also truncation levels in time, namely, $T_{k}=t_{0}(1-2^{-k})$. Integrating this equation in time between $s$ and $t$, where $s\in[T_{k-1},T_{k}]$ and $t\in[T_{k},\infty]$, yields
\[ \int_{M} \theta^{2}_{k}(t)d\textrm{vol}_{g}(x) + 2\int_{s}^{t}\int_{M}\abs{\Lambda^{1/2}\theta_{k}}^{2} \ d\textrm{vol}_{g}(x) \ dt \leq \int_{M} \theta^{2}_{k}(s) \ d\textrm{vol}_{g}(x) .\]

Taking the supremum over $t\geq T_{k}$, 
\[\sup_{t\geq T_{k}}\int_{M}\theta_{k}^{2} \ d\textrm{vol}_{g}(x) + 2 \int_{s}^{\infty}\int_{M} \abs{\Lambda^{1/2}\theta_{k}}^{2} d\textrm{vol}_{g}(x) dt \leq \int_{M} \theta^{2}_{k}(s) \ d\textrm{vol}_{g}(x) .\]
The right hand side dominates the following quantity 
\[E_{k}= \sup_{t\geq T_{k}}\int_{M}\theta_k^{2} \ d\textrm{vol}_{g}(x)+ 2 \int_{T_{k}}^{\infty}\int_{M} \abs{\Lambda^{1/2}\theta_{k}}^{2} d\textrm{vol}_{g}(x) dt.\]
Taking the mean value on the resulting inequality on the interval \\ $s\in[T_{k-1},T_{k}]$ gives
\[ E_{k}\leq \frac{2^{k}}{t_{0}} \int_{T_{k-1}}^{\infty}\int_{M} \theta^{2}_{k} d\textrm{vol}_{g}(x) dt.\]
Notice that for any $x\in M$ such that $\theta_{k}(x)>0$ one also has, by construction, that $\theta_{k-1}(x)\geq 2^{-k}C$. Therefore,
\[ \chi_{\lbrace\theta_{k}>0\rbrace} \leq \left(\frac{2^{k}}{C}\theta_{k-1}\right)^{2/n}.\]
As a consequence of this
\[ \begin{split}
E_{k} &\leq  \frac{2^{k}}{t_{0}} \int_{T_{k-1}}^{\infty} \int_{M} \theta^{2}_{k-1} \chi_{\lbrace\theta_{k}>0\rbrace}  d\textrm{vol}_{g}(x)  dt  \\
 & \leq  \frac{2^{k(1+\frac{2}{n})}}{t_{0}C^{2/n}}\int_{T_{k-1}}^{\infty}\int_{M}\theta^{2(n+1)/n}_{k-1} d\textrm{vol}_{g}(x) dt.
\end{split}\]

Now taking into account that $E_{k-1}$ controls $\theta_{k-1}$ in $L^{\infty}_tL_x^{2}$ and the Sobolev embedding $L^2_tH^{\frac{1}{2}}_x\hookrightarrow L^{2}_tL_x^{2n/(n-1)}$, then H\"older's inequality implies that it also controls $L^{2(n+1)/n}_{t,x}$. Indeed, Sobolev and Poincar\'e inequalities yields
\[\begin{split}
\int_{T_{k-1}}^{\infty}\int_M\theta_{k-1}^{2(n+1)/n}d\textrm{vol}_g(x)dt&\leq\int_{T_{k-1}}^{\infty}\left(\int_M\theta_{k-1}^2d\textrm{vol}_g(x)\right)^{\frac{1}{n}}\left(\int_M\theta_{k-1}^{2n/(n-1)}d\textrm{vol}_g(x)\right)^{(n-1)/n}dt\\
&\leq 2E_{k-1}^{\frac{1}{n}}\int_{T_{k-1}}^{\infty}\int_M|\Lambda^{\frac{1}{2}}\theta_{k-1}|^2d\textrm{vol}_g(x)dt\leq E_{k-1}^{1+1/n}.
\end{split}\]
Therefore, we get the nonlinear recurrence
\[ E_{k} \leq \frac{2^{k(1+	\frac{2}{n})+1}}{t_{0}C^{2/n}} E_{k-1}^{1+\frac{1}{n}},\]
which, for the sake of simplicity, can be rewritten as
\[ E_{k}\leq C' 2^{k(1+2\epsilon)} E_{k-1}^{1+\epsilon},\]
where $C'=C'(M,t_{0},n)$ and $\epsilon=\frac{1}{n}$. We claim that this sequence $E_{k}$ converges to zero if $C$ is large enough (i.e. $C'$ is small enough). Indeed, let us show by induction that $E_k\leq\delta^kE_0$ for $\delta=2^{-(2\epsilon+1)\frac{1}{\epsilon}}< 1$, independent of $k$: we have
\[ \frac{E_{k}}{E_{k-1}} \leq C' 2^{k(1+2\epsilon)} E_{k-1}^{\epsilon} \leq  \delta.\]
For $k=1$, we easily get
\[  C' 2^{(1+2\epsilon)}E_{0}^{\epsilon} \leq \delta\]
choosing the parameter $C'$ sufficiently small. We also have by the induction hypothesis that
\[ E_{k} \leq \delta E_{k-1}\textrm{ and }C'2^{k(1+2\epsilon)}E_{k-1}^{\epsilon}\leq\delta.\]
Therefore
\[\begin{split}
 C' 2^{(k+1)(1+2\epsilon)}E_{k}^{\epsilon} &\leq  C' 2^{(k+1)(1+2\epsilon)}(\delta E_{k-1})^{\epsilon}  \\
 &=C'2^{k(2\epsilon +1)}E_{k-1}^{\epsilon}2^{(2\epsilon +1)}\delta^{\epsilon}\leq \delta.
 \end{split}\]

The non-linear part dissapears since $u$ is divergence free. Furthermore one can mimic the proof for $-\theta$ to achieve the same bound for $\abs{\theta}$. 

\section{\textbf{The barrier functions}}\label{barrier}

The results in this section are instrumental but, nevertheless, they represent the core of the paper. Both lemmas \ref{b1} and \ref{b2} in this section can be interpreted as a quantitative maximum principle for specific boundary elliptic problems at different scales. 

To continue, we need to introduce a piece of notation: in this section we will work on a product space that corresponds to a space variable $x\in M$ times $z\in\mathbb{R}$ (see \S\ref{localenergy} for details) and  $N=n+1$ is its dimension. The arguments in the following sections deal with local properties around some fixed point $x_0\in M$ and a geodesic ball around it $B_g(h)$ of radius $h$ in the metric $g$. The dependence on the point is omitted since our conclusions are uniform due to $M$'s compactness. The usual euclidean metric will be denoted by $g=e$. We will deal with cylinders $(x,z)\in B^*_g(r,h)=B_g(r)\times I(h)$, where $I(h)$ denotes an interval in the variable $z$ of length $h$. Usually, its endpoints will be irrelevant (in the few cases where they are relevant we will point it out explicitly). By a slight abuse of notation we will denote $B^*_g(h,h)$ by $B^*_g(h)$. 

\begin{lemma}\label{b1}
Let the function $b_1$ satisfy
\[\left\{\begin{array}{cl}
(\partial_z^2+\Delta_g)b_1=0&\textrm{in $B^*_g(h)$}\\
b_1=0&\textrm{in $B_g(h)\times\partial I(h)$}\\
b_1=1&\textrm{in $\partial B_g(h)\times I(h)$}
\end{array}\right.\]
Then there exists a $\delta<1$ (independent of the scale $h$) such that for any $x\in B_g(h/2)\times I(h)$
\[b_1(x,z)<\delta+O(h).\]
\end{lemma}

\textsc{Proof of Lemma \ref{b1}:} let $b$ be the euclidean version at scale one
\[\left\{\begin{array}{cl}
(\partial_z^2+\Delta_e)b=0&\textrm{in $B^*_e(1)$}\\
b=0&\textrm{in $B_e(1)\times\partial I(1)$}\\
b=1&\textrm{in $\partial B_e(1)\times I(1)$}
\end{array}\right.\]
Then $b(x/h)$ satisfies the same equation at scale $h$. Define $\delta$ as de supremum of $b(x)$ in $B_e(1/2)\times I(1)$, which is strictly smaller than one by the maximum principle (cf. \cite{LE}, \cite{PW}). We will treat $b_1$ as a perturbation of $b(x/h)$, the difference $u(x)=b_1(x)-b(x/h)$ satisfies
\[\left\{\begin{array}{cl}
(\partial_z^2+\Delta_g)u=O(h^{-1})\frac{\partial b}{\partial\rho}&\textrm{in $B^*_g(h)$}\\
0&\textrm{in $\partial B^*_g(h)$}
\end{array}\right.\]
where $\rho$ is the geodesic radius. One then uses Green's function for the geodesic problem to represent
\[u(x)=O(h^{-1})\int_{B^*_g(h)}G_g(x,y)\frac{\partial b}{\partial\rho}(y/h)d\textrm{vol}_g(y).\]
The integral is bounded (up to a constant dependent on $M$) by
\[\int_{B^*_e(h)}\frac{1}{|x-y|^{N-2}}dy=O(h^2).\]

\begin{figure}[htb]
\centering
\includegraphics[width=90mm]{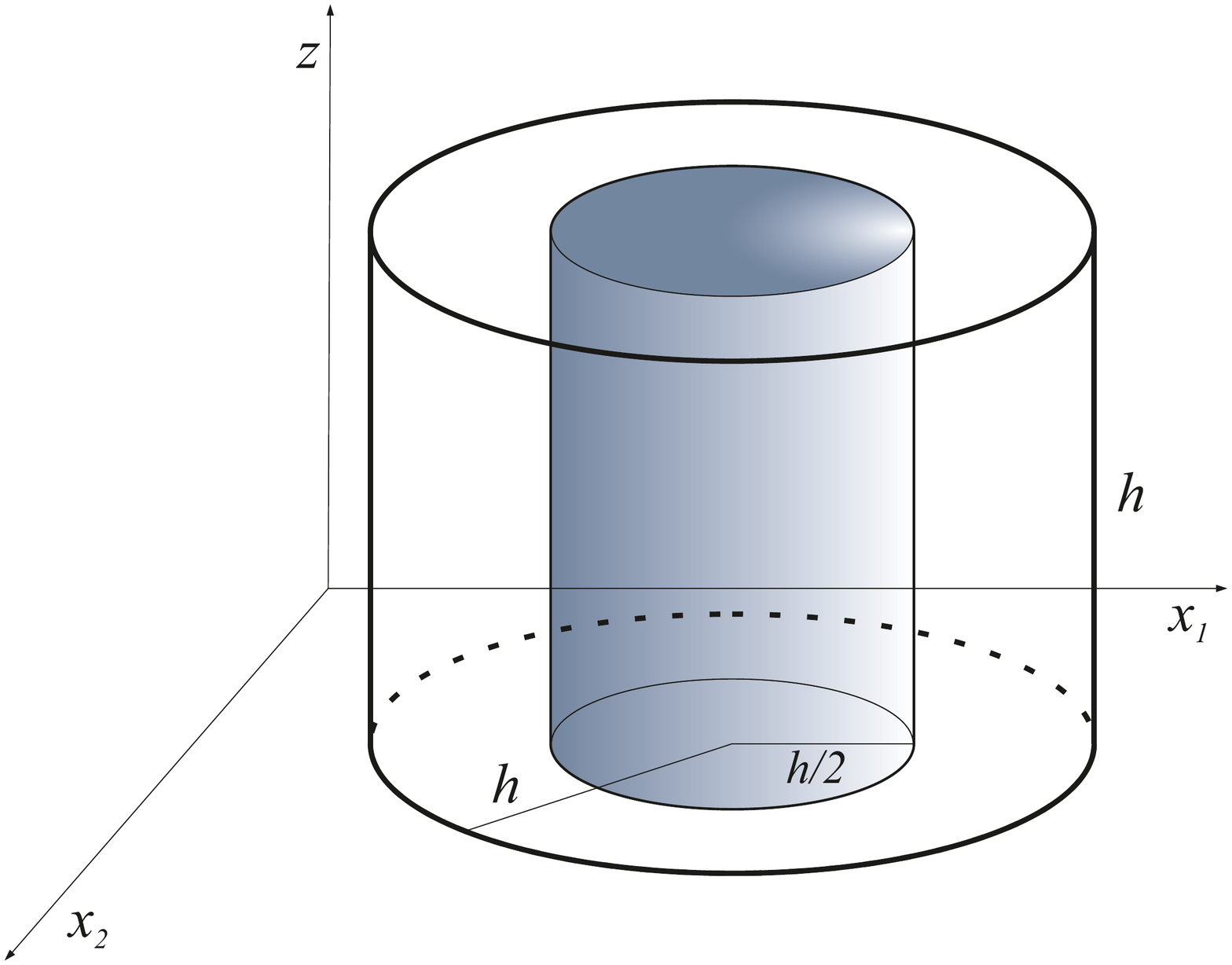}
\caption{Barrier $b_{1}$}\label{fig:b1}
\end{figure}

\textit{Remark:} in the latter bound we used the fact that $G_g(x,y)=O(d(x,y)^{2-N})$ for $N\geq 3$ a fact that follows because the singularity is of that particular order and a maximum principle. The leading term in Hadamard's parametrix shows that the singularity has that prescribed order if $N\geq 3$ (cf. \cite{Ho}). The constants involved depend continuously on the riemannian distorsion of the euclidean metric which can be estimated uniformly due to the assumed compacity. \\

In the following we prove a variant of \ref{b1} that we will explote later, namely: 
\begin{lemma}\label{b2}
Let $h\leq r$. There exist a function $b_2$ such that
\[\left\{\begin{array}{cl}
(\partial_z^2+\Delta_g)b_2=0&\textrm{ in $B^*_g(r,h)$}\\
b_2\geq 0 &\textrm{in $B_g(r)\times\partial I(h)$}\\
b_2=1&\textrm{in $\partial B_g(r)\times I(h)$}
\end{array}\right.\]
and satisfying for $r_1\leq r$
\[\sup_{x\in B_g^*(r_1,h)}b_2(x,z)\leq C\left\{\frac{hr^{N-2}}{(r-r_1)^{N-1}}+r\right\}.\]
\end{lemma}

\begin{figure}[htb]
\centering
\includegraphics[width=90mm]{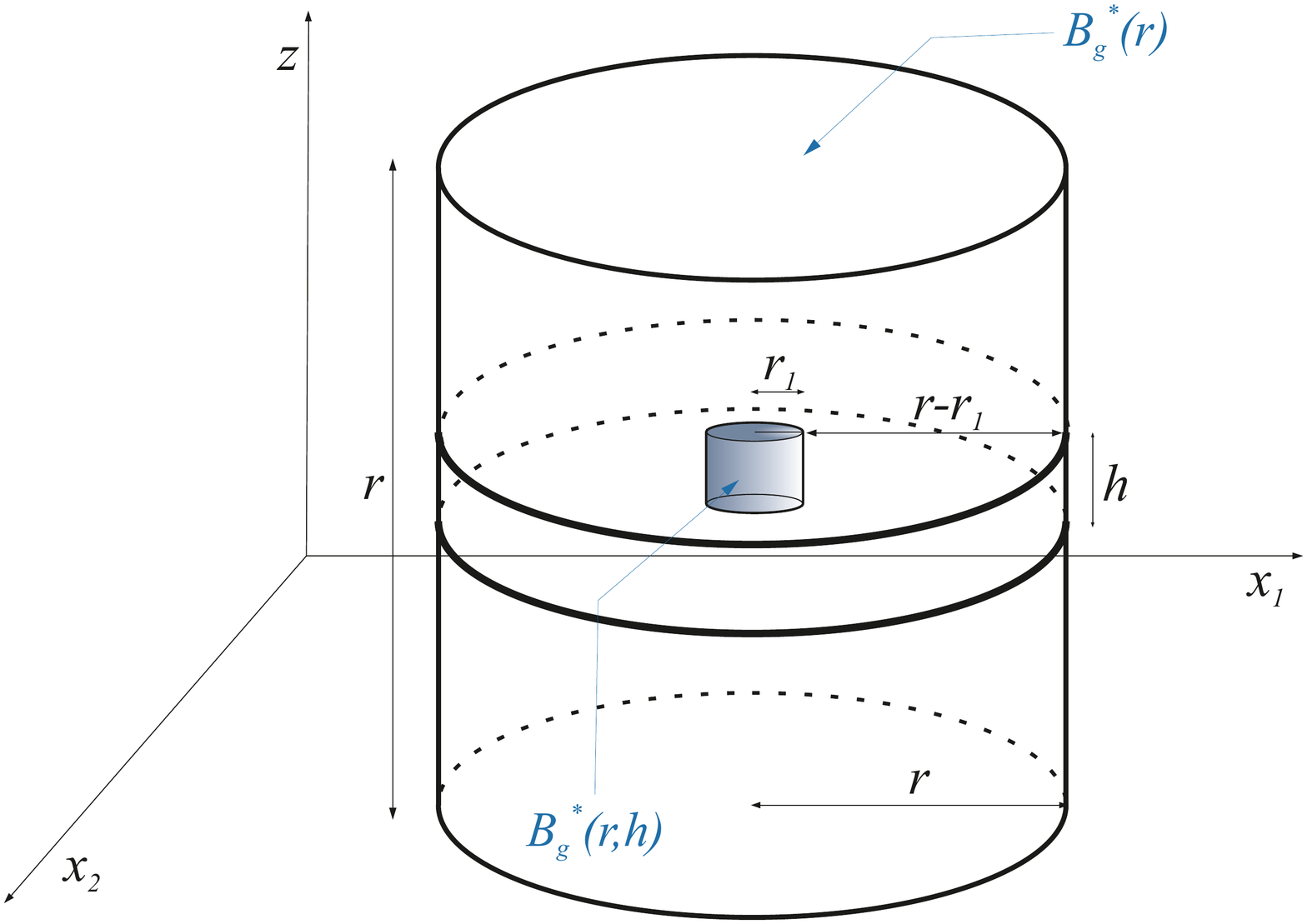}
\caption{Barrier $b_2$}\label{fig:b2}
\end{figure}

\textsc{Proof of Lemma \ref{b2}:} again we will prove it by a perturbation method: first the euclidean and later the general case controlling the difference. If the metric was the euclidean, the use of Green's function estimates yield the result. Indeed, consider $B^*_e(r,h)\subseteq B^*_e(r)$ and let $b$ be the restriction of a function harmonic in $B^*_e(r)$ (we are making the domain larger, see the figure) with non negative boundary values defined to be equal to one near the equator and vanishing outside of it. The maximum principle assures that such a function is non negative and, by construction, satisfies all the assumptions. Finally, observe that integrating against the Poisson kernel $\partial_{\nu}G(x,y)=O(r^{-1}|x-y|^{1-N})$ provides the searched estimate. The need to make the domain larger allows us to rescale the Poisson kernel in the domain $B^*_e(1)$ which implies $\nabla G(x,y)=O(r^{-1}|x-y|^{1-N})$. Otherwise the constant involved might depend on the domain under consideration and, as a matter of fact, on the scale $h$ (cf. Widman \cite{W}, one might in fact round the domain to make it $C^2$ if necessary without affecting the argument above.) Alternatively, one may rescale first to obtain the bound which is invariant upon rescaling and then rescale back (see figure \ref{fig:b2}.)

As in our previous lemma, one treats $u=b_2-b$ as a perturbation that satisfies the following boundary value problem:
\[\left\{\begin{array}{cl}
(\partial_z^2+\Delta_e)u=k\frac{\partial b_2}{\partial\rho}&\textrm{in $B^*_g(r)$}\\
0&\textrm{in the boundary}
\end{array}\right.\]
Here $k$ is a differentiable function independent of $z$ of size $O(r)$ (cf. \cite{Ch2}, theorem 2.17; which can be computed explicitly for the sphere). 

Next we can estimate the derivative using Gauss' divergence theorem. To do that let $x=(\rho,\sigma,z)$ be the cylindrical geodesic coordinates: $\rho$ geodesic distance, $\sigma$ angular direction, $z$ the orthogonal variable. As before one has the following integral representation
\[u(x)=\int_{B^*_e(r)}G_e(x,y)k(y)\frac{\partial b_2(y)}{\partial \rho}dy.\]
Taking into account that $\frac{\partial b_2(y)}{\partial\rho}=\nabla b_2(y)\cdot\sigma$ we can write
\[\begin{split}
u(x)&=\int_{B^*_e(r)}\nabla\cdot\left(G_e(x,y)k(y)b_2(y)\sigma\right)dy-\int_{B^*_e(r)}\nabla G_e(x,y)\cdot\sigma k(y)b_2(y)dy\\
&-\int_{B^*_e(r)}G_e(x,y)\nabla k(y)\cdot\sigma b_2(y) dy-\int_{B^*_e(r)}G_e(x,y)k(y)b_2(y)\nabla\cdot\sigma dy.
\end{split}\]
We may now delete a neighbourhood of $x$ and let it tend to zero to get rid of the singularity of $G_e$ around $x=y$. The first term equals
\[\int_{\partial B^*_e(r)}G_e(x,y)k(y)b_2(y)\sigma\cdot\nu(y) d\sigma(y)\]
which, taking into account that, a priori, $0\leq b_2\leq 1$ by the maximum principle, it can be estimated as $O(r)$. Similarly, the second term is $O(r)$, the third is $O(r^2)$ and the last $O(r)$.

\section{\textbf{Local energy inequality}}\label{localenergy}

In this section we present a local energy inequality (cf. \S\ref{infbound}), that will be used to provide the seeked oscillation decay in \S\S\ref{small}-\ref{big}. Notice that at this stage we know that our weak solution $\theta$ is actually in $L^{\infty}_tL^2_x$ and $L^{\infty}_{t,x}$, therefore by interpolation it follows that $\theta\in L^{\infty}_tL^{2n}_x$, which implies that $u$ is uniformly bounded in $L^{\infty}_tL^{2n}_x$. In the original paper Caffarelli and Vasseur exploited that the drift $u\in L^{\infty}_tBMO_x(\mathbb{R}^n)$, preserved under the natural scaling of the equation. Our approach on the other hand is scale dependent and we will use a localized version instead.

It is useful to think of the fractional Laplace-Beltrami as the boundary value of a derivative through a fractional heat equation, namely
\[\left\{\begin{array}{l}
\partial_zf^*(x,t,z)=-\Lambda^{\alpha}f^*(x,t,z)\\
f^*(x,t,0)=f(x,t)
\end{array}\right.\]
where we denote $z$ this ``time'' variable since we are dealing already with another time variable $t$. Notice $\partial_z f^*(x,t,0)=-\Lambda^{\alpha}f(x,t)$. An additional feature when $\alpha=1$ is that
\[(\partial^2_z+\Delta_g)f^*=0\]
which shows harmonicity for $f^*$, the extension of $f$. This will be a recurrent theme in the sequel. As a consequence of this observation one may use Green's identities in the presence of $\Lambda$, so long as one is willing to work with $f^*$ instead, allowing the treatment of this nonlocal operator as a local one (cf. \cite{CSil, CS}). This idea is exploited deeply in the following Lemma.

\begin{lemma}[Local energy inequality]\label{lemmaenergy}
Let $\theta_k$ satisfy
\[\partial_t\theta_k+u\cdot\nabla_g\theta_k\leq-\Lambda\theta_k\]
and denote $I(z_0)=[0,z_0]$. Let the function $\eta\theta^*_k(x,t,z)$ be vanishing in $M\times[0,\infty)\setminus B_g(h)\times I(z_0)$. Then if $u$ satisfies
\[\displaystyle\sup_{t\in(s,t)}\int_{B_g(h)}|u(x,t)|^{2n}d\textrm{vol}_g(x)\leq Ch^n\]
and $s\leq t$, the following holds
\[\begin{split}
\int_{s}^t\int_{I(z_0)}&\int_{B_g(h)}|\nabla_{x,z}(\eta\theta^*_k)(x,t,z)|^2d\textrm{vol}_g(x)dzdt+\int_{B_g(h)}(\eta\theta_k)^2(x,t)d\textrm{vol}_g(x)\\
&\leq C\left\{\int_{B_g(h)}(\eta\theta_k)^2(x,s)d\textrm{vol}_g(x)+h\int_{s}^t\int_{B_g(h)}|\nabla_{x}\eta\theta_k|^2d\textrm{vol}_g(x)dt\right.\\
&\qquad\quad\quad+\int_{s}^t\int_{I(z_0)}\int_{B_g(h)}|\nabla_{x,z}\eta\theta^*_k|^2d\textrm{vol}_g(x)dzdt\\
&\left.\qquad\qquad\quad+\int_s^t\int_{B_g(h)}(\eta\theta_k)^2(x,t)d\textrm{vol}_g(x)dt\right\}.
\end{split}\]
\end{lemma}

\textit{Remark:} some comments are in order about the notation we have adopted to state the lemma. The function $\theta_k$ in practice will denote a truncation of the weak solution $\theta$ at some level $\ell_k$ (see \S\ref{small} for details). Notice that $\theta_k^*$ refers to the truncation at the same level of $\theta^*$, the extension, which should not be confused with the extension of the truncation $(\theta_k)^*$ (which will never be used in this paper). The gradient $\nabla_{x,z}$ denotes the gradient in the product space $\partial_z+\nabla_g$. We are abusing notation by denoting with $t$ the time variable and the time integration variable, but we hope no confussion arises.

\textsc{Proof of Lemma ~\ref{lemmaenergy}:} as a consequence of subharmonicity of $\theta^*_k$, we get
\[\int_{I(z_0)}\int_{B_g(h)}\eta^2\theta^*_k(\partial_z^2+\Delta_g)\theta_k^*d\textrm{vol}_g(x)dz\geq 0\]
which yields
\[\begin{split}
\int_{I(z_0)}\int_{B_g(h)}|\nabla_{x,z}(\eta\theta^*_k)|^2d\textrm{vol}_g(x)dz&\leq\int_{I(z_0)}\int_{B_g(h)}|\nabla_{x,z}\eta{\theta^*_k}|^2d\textrm{vol}_g(x)dz\\
&\quad+\int_{B_g(h)}\eta^2\theta_k\Lambda\theta_k d\textrm{vol}_g(x).\end{split}\]
After integration by parts in the last integral, only one of the appearing boundary integrals does not vanish. It can be majorized, under the stated hypothesis, by
\[-\frac{1}{2}\left\{\frac{\partial}{\partial t}\int_{B_g(h)}\eta^2\theta_k^2d\textrm{vol}_g(x)+\int_{B_g(h)}\nabla_{x}(\eta^2)\cdot u\theta^2_kd\textrm{vol}_g(x)\right\}.\]
Integrating the resulting equality in the time interval $[s,t]$ one gets
\[\begin{split}
\int_{s}^t\int_{I(z_0)}&\int_{B_g(h)}|\nabla_{x,z}(\eta\theta^*_k)(x,t,z)|^2d\textrm{vol}_g(x)dzdt+\frac{1}{2}\int_{B_g(h)}(\eta\theta_k)^2(x,t)d\textrm{vol}_g(x)\\
&\leq\int_{s}^t\int_{I(z_0)}\int_{B_g(h)}|\nabla_{x,z}\eta\theta^*_k|^2d\textrm{vol}_g(x)dzdt+\frac{1}{2}\int_{B_g(h)}(\eta\theta_k)^2(x,s)d\textrm{vol}_g(x)\\
&\qquad\quad+\left|\int_{s}^t\int_{B_g(h)}\nabla_{x}(\eta^2)\cdot u\theta_k^2 \ d\textrm{vol}_g(x)dt\right|
\end{split}\]
Then the estimate of the last term is the only task left to be done. Using H\"older and Cauchy-Schwarz inequalities we get the bounds
\[\epsilon\int_{s}^t\|\chi_{B_g(h)}\eta\theta_k\|_{L^{\frac{2n}{n-1}}(M)}^2dt+\frac{1}{\epsilon}\int_s^t\|\nabla_{x}\eta\cdot u\theta_k\|_{L^{\frac{2n}{n+1}}(M)}^2dt,\]
where $\epsilon>0$ will be chosen later. Notice that each one of these two terms can be absorbed by some one in the above. Indeed, for the first term we use Sobolev embedding $H^{\frac{1}{2}}\hookrightarrow L^{\frac{2n}{n-1}}$ and the self-adjointness of $\Lambda$ to obtain
\[\begin{split}\int_{s}^t\|\chi_{B_g(h)}\eta\theta_k\|_{L^{\frac{2n}{n-1}}(M)}^2dt&\leq\int_s^t\int_M|\Lambda^{\frac{1}{2}}(\chi_{B_g(h)}\eta\theta_k)|^2d\textrm{vol}_g(x)dt\\ 
&\qquad+\int_s^t\int_{B_g(h)}(\eta\theta_k)^2(x,t)d\textrm{vol}_g(x)dt.\end{split}\]
The second summand is harmless if $\epsilon\leq C$ while the first is bounded by
\[\left.\int_{s}^t\int_{B_g(h)}\eta\theta_k\Lambda(\chi_{B_g(h)}\eta\theta_k)d\textrm{vol}_g(x)dt=-\int_s^t\int_M(\chi_{B_g(h)}\eta\theta_k)^*\partial_z(\chi_{B_g(h)}\eta\theta_k)^*d\textrm{vol}_g(x)dt\right|_{z=0}\]
Now using Green's identities and the decay at infinity, the above integral equals
\[\int_s^t\int_0^{\infty}\int_M|\nabla_{x,z}(\chi_{B_g(h)}\eta\theta_k)^*|^2d\textrm{vol}_g(x)dzdt\]
and Dirichlet principle implies that it is bounded by
\[\int_s^t\int_0^{\infty}\int_M|\nabla_{x,z}(\chi_{B_g(h)}\eta\theta^*_k)|^2d\textrm{vol}_g(x)dzdt\]
Indeed, the harmonic extension is a minimizer for the Dirichlet energy functional and this leads inmediately to
\[\int_{s}^t\|\chi_{B_g(h)}\eta\theta_k\|_{L^{\frac{2n}{n-1}}(M)}^2dt\leq\int_s^t\int_{I(z_0)}\int_{B_g(h)}|\nabla_{x,z}(\eta\theta^*_k)|^2d\textrm{vol}_g(x)dzdt\]
which can be absorbed by the left hand side of the inequality choosing an adequate $\epsilon$. The second term can be handled as follows
\[\begin{split}
\int_s^t\|\nabla_{x}\eta\cdot u\theta_k\|_{L^{\frac{2n}{n+1}}(M)}^2dt&=\int_s^t\left(\int_M|\nabla_{x}\eta\cdot u\theta_k|^{\frac{2n}{n+1}}d\textrm{vol}_g(x)\right)^{\frac{n+1}{n}}dt\\
&\leq\int_s^t\left(\int_{B_g(h)}|u|^{2n}d\textrm{vol}_g(x)\right)^{\frac{1}{n}}\left(\int_{B_g(h)}|\nabla_{x}\eta\theta_k|^2d\textrm{vol}_g(x)\right)dt\\
&\leq \|u\|_{L^{\infty}_tL^{2n}_x(B_g(h))}^2\int_s^t\int_{B_g(h)}|\nabla_{x}\eta\theta_k|^2d\textrm{vol}_g(x)dt
\end{split}\]
where we have used H\"older's inequality and the fact that $\nabla\eta$ is supported in $B_g(h)$.

\section{\textbf{H\"older regularity}}\label{holder}

This section will deal with the H\"older continuity of weak solutions. The approach is based on the decrease of the $L^{\infty}$ norm of either the positive part or the negative part of $\theta$, which implies a decrease in the oscillation. The proof is subdivided in two stages each one containing a step towards the result. In the first we study the decrease under small mean energy hypothesis while in the second we remove such a restriction. This is reminiscent of De Giorgi's work on Hilbert's 19th problem.

\subsection{Small mean energy}\label{small}

To state precisely the concrete piece of the proof that we will be dealing with in this section, it is convenient to introduce the notation $Q_g(h)$ to denote the pairs $(x,t)$ such that $x\in B_g(h)$ and $t\in t^*+I(h)$, where $t^*\geq t_0$. Following \S\ref{infbound} we use $Q_g^*(h)$ to denote the set of $(x,t,z)$ where $(x,t)\in Q_g(h)$ and $z\in I(h)$. Notice that we are not accurate about the precise position of the time interval, but we only care about its length. In the sequel time intervals will be chosen carefully.

\begin{proposition}\label{propsmall}
For $h$ small enough, there exist $\epsilon>0$ and $\gamma<1$ (both independent of the scale $h$) so that for any solution $\theta$ satisfying
\[\int_{Q_g(2h)}(\theta)_+^2d\textrm{vol}_g(x)dt\leq\epsilon \int_{Q_g(2h)}d\textrm{vol}_g(x)dt\]
and
\[\int_{Q^*_g(2h)}(\theta^*)_+^2d\textrm{vol}_g(x)dt\leq\epsilon \int_{Q^*_g(2h)}d\textrm{vol}_g(x)dzdt\]
one also has 
\[\|\theta_+\|_{L^{\infty}(Q_g(h))}\leq\gamma\|\theta^*\|_{L^{\infty}(Q^*_g(2h))}.\]

\end{proposition}

Some comments are needed before proceeding to the proof itself. The statement is written in terms of $L^{\infty}$ bounds because they are related to the oscillation decrease as follows: instead of $\theta^*$ one may consider $\theta^*-a$ for any arbitrary constant $a$, the resulting function has the same oscillation and satisfies the same drift equation (it is a feature of the proof of Theorem \ref{II} that the active scalar dependence is quite irrelevant). As a consequence of this one may choose $a$ in such a way that the $L^{\infty}$-norm of $(\theta^*-a)_+$ and the oscillation of $\theta^*-a$ are comparable and the decrease on the oscillation is strictly smaller than one. In fact, one may choose $a$ so that $\|(\theta^*-a)_+\|_{L^{\infty}(Q^*_g(h))}$ equals $\|(\theta^*-a)_-\|_{L^{\infty}(Q^*_g(h))}$ and  hence both are precisely $\frac{1}{2}\textrm{osc}_{Q^*_g(h)}\theta^*$, where $(f)_-=-(-f)_+$. Provided the conclusion of the proposition is true one might bound
\[\begin{split}
\textrm{osc}_{Q_g(h/2)}\theta&=\textrm{osc}_{Q_g(h/2)}(\theta-a)\leq\|(\theta-a)_+\|_{L^{\infty}(Q_g(h/2))}+\|(\theta-a)_-\|_{L^{\infty}(Q_g(h/2))}\\
&\leq\gamma\|(\theta^*-a)_+\|_{L^{\infty}(Q_g(h)}+\|(\theta-a)_-\|_{L^{\infty}(Q_g(h))}\\
&\leq\frac{1}{2}\left(1+\gamma\right)\textrm{osc}_{Q^*_g(h)}(\theta^*-a)=\frac{1}{2}\left(1+\gamma\right)\textrm{osc}_{Q^*_g(h)}\theta^*.
\end{split}\]
Hence, we have shown as a byproduct that the oscillation would decrease by $(1+\gamma)/2<1$.

We will proceed by a rather tricky induction process involving some local energy quantities in the same spirit as in De Giorgi's tecnique, but needing to control several boundary terms due to the nonlocality. Since the proof is quite technical and intrincated, let us expose first the general plan of how to achieve the nonlinear inequality for the local energy $E_{k}$  we are aiming to. In order to clarify the exposition, we state certain claims whose proof will be postponed to the end of this digression.

Fix some $\gamma<1$ to be specified later. Let us denote by $\theta_k$ and $\theta_k^*$ the positive part of the trucations of $\theta$ and $\theta^*$ respectively at the level
\[\ell_k=\|\theta^*\|_{L^{\infty}(Q^*_g(2h))}\left(1-(1-\gamma)\frac{{1+2^{-k}}}{2}\right)\] 
Let $\eta_k$ be a smooth bump function supported in $B_g(h(1+2^{-k}))$ identically one in $B_g(h(1+2^{-k-1}))$. We will find a nonlinear inequality for the local energy
\[\begin{split}
E_k&=\sup_{t\in t^*+[-h2^{-k-1},h]}\int_M(\eta_k\theta_k)^2(x,t)d\textrm{vol}_g(x)\\
&\qquad\qquad+\int_{t^*+[-h2^{-k-1},h]}\int_{I(h\delta^k)}\int_M|\nabla_{x,z}(\eta_k\theta_k^*)(x,t,z)|^2d\textrm{vol}_g(x)dzdt
\end{split}\]
where we assumed $t^*-t_0\leq h/4$ (cf. \cite{CV}). Otherwise one may shrink the intervals in both definitions, using for the lower extreme $t^*-h2^{-k-k_0}$ instead, for some appropiate $k_0\geq 0$. The constant $\delta<1$ is a small parameter to be selected later independently of $h$. This choice for $E_k$ is motivated by the proof \S\ref{infbound} and the local energy inequality from \S\ref{localenergy} (cf. inequality $(*)$ below). Notice that $E_k$ decrease. Recall that we can not afford fixing a reference scale, and use the iterative scaling arguments to deal with other finer scales. Therefore we need to keep track upon the scale during the proof, and this fact will be crucial for our argument to finally work properly.

Taking mean value for $s\in t^*-[h2^{-k-1},h2^{-k}]$ in the local energy inequality of lemma \ref{lemmaenergy}) for $z_0=h\delta^{k}$, $\eta=\eta_k$ and any $t\in t^*+[-h2^{-k-1},h]$ one obtains
\begin{equation*}\label{localmedia}\begin{split}
&\int_{B_g(2h)}(\eta_k\theta_k)^2(x,t)d\textrm{vol}_g(x)\qquad\quad\quad\qquad\quad\qquad\qquad\qquad\qquad\quad\qquad\quad\qquad(*)\\
&\qquad+\int_{t^*+[-h2^{-k-1}, t]}\int_{I(h\delta^k)}\int_{B_g(2h)}|\nabla_{x,z}(\eta_k\theta^*_k)|^2d\textrm{vol}_g(x)dzdt\\
&\qquad\qquad\leq C\left\{\frac{2^k}{h}\int_{t^*+[-h2^{-k},h]}\int_{B_g(2h)}(\eta_k\theta_k)^2(x,s)d\textrm{vol}_g(x)ds\right.\\
&\qquad\qquad\qquad\quad+h\int_{t^*+[-h2^{-k},h]}\int_{B_g(2h)}|\nabla_{x}\eta_k\theta_k|^2d\textrm{vol}_g(x)dt\\
&\qquad\qquad\qquad\qquad\quad+\left.\int_{t^*+[-h2^{-k},h]}\int_{I(h\delta^k)}\int_{B_g(2h)}|\nabla_{x,z}\eta_k\theta^*_k|^2d\textrm{vol}_g(x)dzdt\right\}.
\end{split}
\end{equation*}
The constant $C$ does not depend on the scale $h$ nor on the truncation step $k$ (cf. recall remark ~\ref{constantes} about the use of constants). Taking supremum on $t$ and using the elementary bound $|\nabla\eta_k|\leq C\frac{2^k}{h}\eta_{k-1}$, which can be assumed to be true for a certain construction of the bump functions, one obtains
\[\begin{split}
E_k&\leq C\frac{2^{2k}}{h^2}\left\{h\int_{t^*+[-h2^{-k},h]}\int_{B_g(2h)}(\eta_{k-1}\theta_k)^2(x,s)d\textrm{vol}_g(x)ds\right.\\
&\left.\qquad\quad\quad\quad\quad+\int_{t^*+[-h2^{-k},h]}\int_{I(h\delta^k)}\int_{B_g(2h)}(\eta_{k-1}\theta^*_k)^2d\textrm{vol}_g(x)dzdt\right\}.
\end{split}\]
By construction for any $x$ such that $\theta_k(x)>0$, one has $\theta_{k-1}(x)\geq C(1-\gamma)2^{-k-2}$. Using $\chi_{\{\eta_{k-1}>0\}}\chi_{\{\theta_k>0\}}\leq C\frac{2^{k}}{1-\gamma}\theta_{k-1}\eta_{k-2}$ in the above we get the bound
\[\begin{split}
E_k&\leq C\frac{2^{3k}}{h^2(1-\gamma)^{2/n}}\left\{h\int_{t^*+[-h2^{-k},h]}\int_M(\eta_{k-2}\theta_{k-1})^{2\frac{n+1}{n}}d\textrm{vol}_g(x)dt\right.\\
&\left.\qquad\qquad\qquad\qquad\qquad+\int_{t^*+[-h2^{-k},h]}\int_{I(h\delta^k)}\int_M(\eta_{k-2}\theta_{k-1}^*)^{2\frac{n+1}{n}}d\textrm{vol}_g(x)dzdt\right\}.
\end{split}\]
This will be estimated from above in the same nonlinear way as was done in \S\ref{infbound}. To do so let us first claim that $\eta_{k-2}\theta_{k-2}^*=0$ provided $z\in [h\delta^{k-1},h\delta^{k-2}]$ (cf. Lemma \ref{conditions})
\[\begin{split}
\int_{I(h\delta^{k-2})}\int_M|\nabla_{x,z}(\eta_{k-2}\theta_{k-2}^*)|^2d\textrm{vol}_g(x)dz&=\int_0^{\infty}\int_M|\nabla_{x,z}(\eta_{k-2}\theta^*_{k-2}\chi_{I(h\delta^{k-2})})|^2d\textrm{vol}_g(x)dz\\
&\geq\int_0^{\infty}\int_M|\nabla_{x,z}(\eta_{k-2}\theta_{k-2})^*|^2d\textrm{vol}_g(x)dz\\
&=\left.-\int_M\eta_{k-2}\theta^*_{k-2}\partial_z(\eta_{k-2}\theta^*_{k-2})d\textrm{vol}_g(x)\right|_{z=0}\\
&=\int_M|\Lambda^{\frac{1}{2}}(\eta_{k-2}\theta_{k-2})|^2d\textrm{vol}_g(x)
\end{split}\]
where we have used the aforementioned claim, decay at infinity, harmonicity and Green's identities. This shows that $E_{k-2}$ dominates the norms  $L^2H^{\frac{1}{2}}\hookrightarrow L^2L^{2n/(n-1)}$ and $L^{\infty}L^2$, using H\"older's inequality as in \S\ref{infbound} we obtain that
\[\int_{t^*+[-h2^{-k},h]}\int_M(\eta_{k-2}\theta_{k-2})^{2(n+1)/n}d\textrm{vol}_g(x)dt\leq E_{k-2}^{1+1/n}\]
This suggests that an inequality of the type
\[E_{k}\leq C\frac{2^{3k}}{(1-\gamma)^{2/n}h}E_{k-2}^{1+1/n}\]
might hold true, which is almost the kind of non linear inequality we would like to use.\footnote{Notice the $h$ in the denominator is harmless since we are working with the hypothesis in the mean. This is a reminiscence of the standard dimensional analysis of physical quantities.} The estimate of the remaining term, $\|\eta_{k-2}\theta^*_{k-2}\|^{2(n+1)/n}_{L^{2(n+1)/n}(M)}$, can be reduced to the above. However, since this is not inmediate, let us show first that for any $t$
\[\theta_{k+1}^*(x,t,z)\leq(\eta_k\theta_k)^*(x,t,z)\textrm{ for any $(x,z)\in B^*_g(h(1+2^{-k-1}), h\delta^k)$}\]
holds provided the claim is true. Indeed, by harmonicity one has the bound
\[\theta_k^*(x,t,z)\leq\int_M\eta_k(y)\theta_k(y,t)G(x,y,z)d\textrm{vol}_g(y)+\|\theta^*_k\|_{L^{\infty}(Q^*_g(2h))}b_2(x)\]
for any pair $(x,z)\in B^*_g(h(1+2^{-k-1}),h\delta^k)$, which follows using the maximum principle in the cylinder to majorize $\theta^*_k$ by $(\eta_k\theta_k)^*$ in the bottom part of the cylinder (i.e. for $x\in B_g(h(1+2^{-k-1}))$ and $z=0$). On the other hand, the claim allows to disregard the upper part which is bounded by zero. In the rest of the boundary we use the barriers introduced in \S\ref{barrier}, which by construction verify (cf. Lemma \ref{b2})
\[|b_2(x)|\leq C\left\{\delta^{k}(1+2^{-k})^{n-2}2^{k(n-1)}+h(1+2^{-k})\right\}\]
which is smaller than $(1-\gamma)2^{-k-2}$ provided $h$ and $\delta$ are small enough independently of $k$. From now on we will suppose that $h$ and $\delta$ are such that the above holds true. To continue let us observe that
\[\theta^*-\ell_{k+1}=\theta^*-\ell_k-(1-\gamma)2^{-k-1}\|\theta^*\|_{L^{\infty}(Q^*_g(2h))}\]
 from which one gets the inequality
\[\theta_{k+1}^*(x,t,z)\leq\int_M\eta_k(y)\theta_k(y,t)G(x,y,z)d\textrm{vol}_g(y)\textrm{ for }x\in B_g(h(1+2^{-k-1}))\]
which yields $\eta_{k+1}\theta^*_{k+1}\leq(\eta_k\theta_k)^*$. Using this fact we get the estimate
\[\begin{split}
\int_{t^*+[-h2^{-k},h]}&\int_{I(h\delta^k)}\int_M(\eta_{k-2}\theta_{k-1}^*)^{2\frac{n+1}{n}}d\textrm{vol}_g(y)dzdt\\
&\leq \int_{t^*+[-h2^{-k},h]}\int_{I(h\delta^k)}\int_M|(\eta_{k-3}\theta_{k-3})^*|^{2\frac{n+1}{n}}d\textrm{vol}_g(y)dzdt\\
&\qquad\leq \int_{t^*+[-h2^{-k},h]}\int_{I(h\delta^k)}\int_M(\eta_{k-3}\theta_{k-3})^{2\frac{n+1}{n}}d\textrm{vol}_g(y)dzdt.
\end{split} \]
where we have applied Jensen's inequality and the identity $\int_MG(x,y,z)d\textrm{vol}_g(y)=1$. This is already known to be bounded nicely in terms of $E_k$ as above, provided the claim holds. We conclude the digression by observing that as a consequence of the decreasing character of the energies one would get the following nonlinear inequality
\begin{equation}\label{plausible}
E_{k}\leq C\frac{2^{3k}}{(1-\gamma)^{2/n}h}E_{k-3}^{1+1/n}
\end{equation}

The next step is to prove our claim provided some extra hypothesis is fulfilled. This will be helpful to close the induction later on.

\begin{lemma}\label{conditions}
For $k\geq 0$, the following statement holds:
\begin{equation}\label{set}
 \theta_{k+1}^*(x,t,z)= 0\textrm{ for any $(x,z)\in B_g(h(1+2^{-k}))\times(h\delta^{k+1},z_0)$}
\end{equation}
provided that $\theta^*_{k}(x,t,z_0)=0$ and the energy $E_k$ satisfies 
\begin{equation}\label{smallness}
 E_kh^{-n}\delta^{-2n(k+1)}\leq C^2(1-\gamma)^22^{-2(k+1)}
\end{equation}
\end{lemma}

\textsc{Proof of lemma ~\ref{conditions}:} we would like to bound $(\eta_k\theta_k)^*$ in the preceding discussion  by $C(1-\gamma)2^{-k-2}$ which, intertwined with the arguments above, will be enough for our purposes. Now if $(x,z)\in B_g(h(1+2^{-k}))\times (h\delta^{k+1}, h\delta^{k})$ and $t\in t^*+I(h)$ then we can estimate the other term appropiately, indeed:
\[\int_M\eta_k(y)\theta_k(y,t)G(x,y,z)d\textrm{vol}_g(x)\leq\sqrt{E_k}\|G(x,y,z)\|_{L^2_x(M)}.\]
In the next, we will make use of the following expansion of the fractional heat kernel (cf. \cite{CM})
\[G(x,y,z)=\sum_{i=0}^{\infty}e^{-\lambda_iz}Y_i(x)Y_i(y)\]
which combined with the local Weyl estimates (cf. Theorem 3.3.1, \cite{So} p. 53)
\[\sum_{\lambda_i\leq\lambda}Y_i(x)Y_i(y)=O(\lambda^{n})\]
and a straightforward summation by parts yields
\[G(x,y,h\delta^{k+1})\leq C(h\delta^{k+1})^{-n}\]
from which $\|G(x,y,h\delta^{k+1})\|^2_{L^2_x(M)}\leq C\delta^{-2n(k+1)}h^{-n}$. As a consequence of the hypothesis $(\eta_k\theta_k)^*$ is smaller than $C(1-\gamma)2^{-k-2}$. This altogether with 
\[\theta_{k+1}^*\leq\theta_k^*-(1-\gamma)2^{-k-1}\|\theta^*\|_{L^{\infty}(Q^*_g(2h))}\]
shows that $\theta_{k+1}^*$ can not be positive in $B_g(h(1+2^{-k}))\times (h\delta^{k+1}, z_0)$. \\ 

\textit{Remark:} we need to choose $z_{0}$ in a specific way, say $z_{0}= h/2$, in order to start the inductive procedure (indicated below in the proof of Proposition \ref{propsmall}). Once we have used this, the parameter $z_{0}$ will have the form $\delta^{k}h$ in the $k$th step of the induction process (we hope that it does not create any confusion.)

We believe that Lemma ~\ref{conditions} is a convenient intermediate step to write the induction neatly. Next we proceed to the uncover the details of the proof.

\textsc{Proof of proposition \ref{propsmall}:}  we will show that if $\epsilon$ is small enough one can choose a positive $\beta<1$, independent of $h$, such that 
\[ E_k\leq \beta^{k}h^n \text{ \ holds for any \ \ } k\geq 0. \]
In particular $E_{k}$ tends to zero, proving the statement. The geometrical decay of this ansatz is very convenient in order to check that the hypothesis (\ref{smallness}) imposed to $E_{k}$ of Lemma ~\ref{conditions} is satified (recall that this kind of behaviour is quite plausible for sequences satisfying a non linear inequality (\ref{plausible}) (cf. \S\ref{infbound})). 

To that end we choose some $\beta<1$ satisfying the following smallness condition,
\[\beta^{k}\leq \delta^{2n(k+1)} C(1-\gamma)^24^{-k-1}\]
provided $\delta$ and $\gamma$ are fixed already. Moreover we will also impose
\[C\frac{2^{3k}}{(1-\gamma)^{2/n}}\beta^{\frac{k-3}{n}-4}<1\]
which is only useful when $k\geq 4n+4$. 

Next we will prove by induction the following predicate:
\[P(k): E_{k}\leq\beta^{k}h^n\textrm{ and }\eta_k\theta^*_k=0\textrm{ in the set (\ref{set}) of Lemma \ref{conditions}}.\]
Due to the shift in the nonlinear inequality (\ref{plausible}) and our previous arguments, we already know that if predicates $P(k-3)$, $P(k-2)$, $P(k-1)$ are satisfied then $P(k)$ is also fulfilled, provided $k\geq 4n+4$. Therefore our work is reduced to check that $P(0)$, $P(1)$,\ldots, $P(4n+3)$ are satisfied. The first part of the predicate is quite straightforward. Indeed, using the local energy inequality derived in \S\ref{localenergy}, we can take $\epsilon$ verifying
\[ \epsilon < C \beta^{4n+3} \] 
and hence $E_{k}\leq \beta^{4n+3}h^{n}$, for $0\leq k \leq 4n+3$, as required. 

Next, let us realize that if we prove the second part of the statement for $k=0$, we would be done. Indeed, appealing to Lemma \ref{conditions}, we would prove $P(1)$. Afterwards using $P(1)$ and smallness on $E_1$, we deduce $P(2)$. Similarly one gets $P(3)$. Therefore the induction process works nicely without any further assumptions using the non linear inequality (\ref{plausible}) provided $\beta$ is smaller than a threshold quantity, which is independent of the scale $h$ as we specified previously.

\begin{figure}[H]
\centering
\includegraphics[width=90mm]{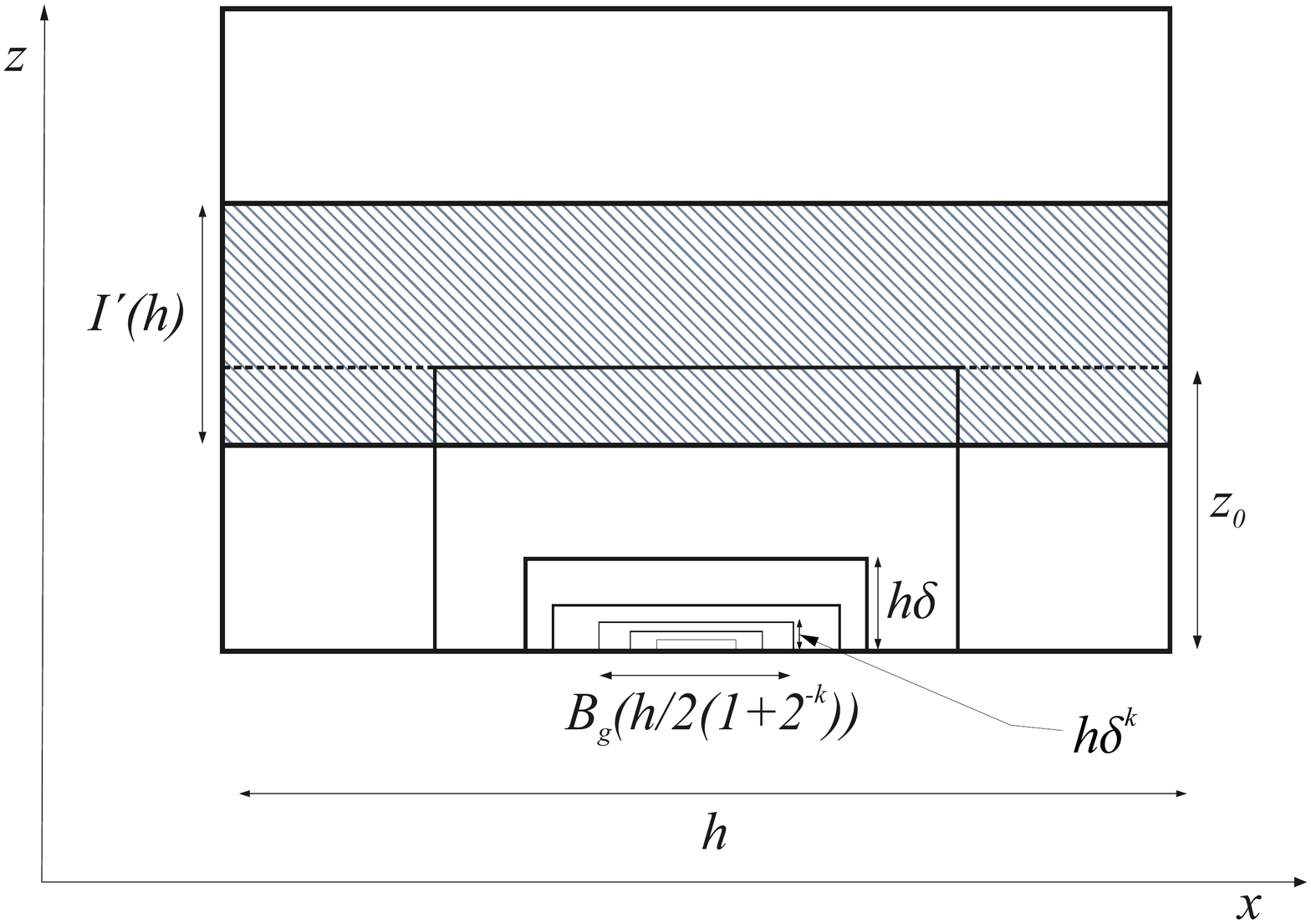}
\caption{Initial step and iterative procedure.}\label{fig:b3}
\end{figure}

Let us deal now with that initial codition. The maximum principle allows us to bound $\theta^*$ in $B_g(h)\times I(h)$ as follows
\[\begin{split}
\theta^*(x,t,z)&\leq\int_M\theta(y,t)\chi_{B_g(h)}(y)G(x,y,z)d\textrm{vol}_g(y)\\
&\qquad\qquad+\|\theta^*\|_{L^{\infty}(Q^*_g(2h))}\left(b_1(x,z)+\frac{z}{h}\right)
\end{split}\]
where the barrier $b_1$ has been constructed in \S\ref{b1}. The only problematic term is the first one since the second term can be handled using Lemma ~\ref{b1} and the third term can be bounded easily due to its linearity. By H\"older's inequality we get
\[\int_M\theta(y,t)\chi_{B_g(h)}(y)G(x,y,z)d\textrm{vol}_g(y) \leq C \norm{\theta^*}_{L^{\infty}(Q_g^{*}(2h))} \norm{G(x,y,z)}_{2}\norm{\chi_{B_{g}(h)}}_{2}\]
which is less than $C(M)2^n$. Using Weyl's law asymptotics, positivity of the characteristic function and fractional heat kernel estimates one gets
\[|\theta^*(x,t,z)| \leq\gamma^*\|\theta^*\|_{L^{\infty}(Q^*_g(2h))} \]
where $x\in B_{g}(h)$, $z\in I'(h)$ and $\gamma^*<1.$ This shows that the hypothesis of Lemma \ref{conditions} are satisfied for $\theta_{0}^*$ with $\gamma=\frac{1}{3}(4\gamma^*-1)$ and $z_0=\frac{h'}{2}$, where we can take $z_0\in I'(h)=[\frac{h}{3},\frac{2h}{3}]$ (see figure \ref{fig:b3}).

\textit{Remark:} in the estimate above we are assuming that $C(M)2^n$ is small enough. But if that is not the case, one may use instead of balls decreasing by half, balls decreasing by some fixed small quantity $c$. Then, the estimate above will have the form $C(M)2^nc^{n/2}$ and we can choose $c$ so that the above estimate is indeed strictly smaller than one. And that is all we need.

\subsection{Arbitrary energy}\label{big}

The purpose of this section is to free Proposition \ref{propsmall} from its small mean energy requirement. To do so we prove a version of De Giorgi's isoperimetric inequality following closely the argument in \cite{CV} though it needs careful adaptation to avoid problems with the different scales. Let us first introduce some convenient notation: denote by $Q_g^*(h)$ the cube of all $(x,t,z)\in B_g(h)\times I(h)\times I(h)$, $Q_g(h)$ the set $(x,t)\in B_g(h)\times I(h)$ and $|A|$ the measure of $A$ in the product. It would be useful while reading this section to keep in mind that $B_g(h)$, $B_g^*(h)$, $Q_g^*(h)$ are approximately of order $h^n$, $h^{n+1}$, $h^{n+2}$, respectively.

De Giorgi's isoperimetric inequality, in our setting, will relate the measures of the following sets
\[\left\{\begin{array}{l}
\mathcal{A}(t)=\{(x,z)\in B_g^*(h):\theta^*(x,t,z)\leq 0\}\\
\mathcal{B}(t)=\left\{(x,z)\in B_g^*(h):\theta^*(x,t,z)\geq \frac{1}{2}\|\theta^*\|_{L^{\infty}(Q_g^*(h))}\right\}\\
\mathcal{C}(t)=\left\{(x,z)\in B_g^*(h):0<\theta^*(x,t,z)< \frac{1}{2}\|\theta^*\|_{L^{\infty}(Q_g^*(h))}\right\}
\end{array}\right.\]
It reads as follows:

\begin{lemma}[De Giorgi's isoperimetric inequality] \label{degiorgi}
For $h$ small enough, the following inequality holds for any function $\theta^*\in H^1(B_g(h))$
\begin{equation}\label{isoperimetrica}
|\mathcal{A}(t)||\mathcal{B}(t)|\leq C |\mathcal{C}(t)|^{\frac{1}{2}}K^{\frac{1}{2}}h^{n+2}
\end{equation}
where $K=\|\nabla_{x,z}\theta^*\|^2_{L^2(B^*_g(h))}$. 
\end{lemma}

\textsc{Proof of Lemma ~\ref{degiorgi}:} proceed as in \cite {CV}, \cite{CV2} being careful to keep the small scale dependence. 

Now we can state the main result:

\begin{lemma}\label{isononlocal}
For any $\epsilon>0$ small enough\footnote{Depending on $M$, $\theta_0$, $t_0$.} (independent of $h$) there exists  $\delta=\delta(\epsilon)>0$ such that for any weak solution $\theta$ satisfying
\[|\{(x,t,z)\in Q^*_g(2h):\theta^*(x,z,t)\leq 0\}|\geq \frac{1}{2}|Q^*_g(2h)|.\]
We have that the hypothesis
\[\left|\left\{(x,t,z)\in Q^*_g(2h):0<\theta^*<\frac{1}{2}\|\theta^*\|_{L^{\infty}(Q^*_g(h))}\right\}\right|\leq\delta|Q^*_g(2h)|\]
implies
\[\int_{Q_g(h)}\left(\theta-\frac{1}{2}\|\theta^*\|_{L^{\infty}}\right)_+^2d\textrm{vol}_g(x)dt\leq \epsilon|Q_g(h)|\]
and
\[\int_{Q^*_g(h)}\left(\theta^*-\frac{1}{2}\|\theta^*\|_{L^{\infty}}\right)_+^2d\textrm{vol}_g(x)dt\leq \epsilon|Q^*_g(h)|.\]
\end{lemma}

Before proceeding to the proof itself let us glimse the rough idea behind it: let us suppose that both $\mathcal{B}(t)$ and $\mathcal{C}(t)$ are smaller than $\gamma h^{n+1}$
\[\int_{B^*_g(h)}(\theta^*)_+^2(t)dxdz\leq\|\theta^*\|_{L^{\infty}(B^*_g(h))}^2(|\mathcal{B}(t)|+|\mathcal{C}(t)|)\leq\gamma Ch^{n+1}\]
(cf. the Remark concerning the use of constants in \S\ref{infbound}.) Integrating the identity
\begin{equation}\label{isoperometrica2}
\int_{B_g(h)}\theta_+^2d\textrm{vol}_g(x)=\int_{B_g(h)}(\theta^*)_+^2d\textrm{vol}_g(x)-2\int_0^z\int_{B_g(h)}\theta_+^*\partial_z\theta^*d\textrm{vol}_g(x)d\bar{z}
\end{equation}
 for $z\in I(h)$ one gets
\[h\int_{B_g(h)}\theta_+^2(x,t)dx\leq\int_{B_g^*(h)}\theta_+^*(x,t,z)^2dxdz+\int_0^h\int_0^z\int_{B^*_g(h)}\theta_+^*(t)\partial_z\theta^*dxd\bar{z}dz.\]
The first term can be bounded by $O(\gamma h^{n+1})$, using the previous inequality; the second, applying Fubini and Cauchy-Schwarz is bounded by 
\[h\|\theta_+^*(t)\|_{L^2(B_g^*(h))}\|\partial_z\theta^*\|_{L^2(B^*_g(h))}\leq\gamma h^{n+1}.\]
 Notice that the $L^2$-gradient norm might be expected to be of size $h^{(n-1)/2}$, from dimensional considerations. Summarizing: if the above argument works we may  integrate the resulting inequality for all times $t\in I(h)$ achieving
\[\int_{Q^*_g(h)}\theta_+^2(x,t) d\textrm{vol}_g(x)dt=O(\sqrt{\gamma} h^{n+2}).\]
Notice that this estimate is stronger than the one we intend to proof. In fact, the assumption above should not be expected to hold for any $t\in I(h)$, the proof will show that an elaboration of the aforementioned argument reaching control on the size of the corresponding sets actually holds for most of the times $t\in I(h)$. The remaining times, for which it does not hold, will have a controlled size. 

\textsc{Proof of Lemma ~\ref{isononlocal}:} from De Giorgi's isoperimetric inequality (\ref{isoperimetrica}) one notice that one may control $\mathcal{B}(t)$'s smallness if one knows for an appropiate $K$ that $\mathcal{A}(t)$ is big, while $\mathcal{C}(t)$ is small (by hypothesis). Let us introduce a subset of times for which we do expect some control of $\mathcal{B}$ due to De Giorgi's isoperimetric inequality provided we manage to prove $\mathcal{A}(t)$ is big enough:
\[\mathcal{T}=\left\{t\in I(h):\int_{B_g^*(h)}|\nabla\theta^*_+|^2d\textrm{vol}_g(x)dz\leq K\textrm{ and } |\mathcal{C}(t)|^{\frac{1}{2}}\leq2\epsilon^3h^{\frac{n+1}{2}}\right\}\]
The complement of this set is small in $I(h)$, in the sense that it is smaller than $\epsilon h/2$, choosing
\[K=\frac{4}{\epsilon h}\int_{Q^*_g(h)}|\nabla_{x,z}\theta^*_+|^2d\textrm{vol}_g(x)dzdt.\]
Indeed, define $\delta=\epsilon^8$, one may obtain the following weak bound 
\[\left|\left\{t\in I(h):|\mathcal{C}(t)|^\frac{1}{2}\geq2\epsilon^3h^{\frac{n+1}{2}}\right\}\right|\leq\frac{1}{4\epsilon^6h^{n+1}}\int_{I(h)}|\mathcal{C}(t)|dt\leq\epsilon^2h/4\]
Furthermore, the control on the remaining condition is provided by the following weak bound 
\[\left|\left\{t\in I(h):\int_{B_g^*(h)}|\nabla_{x,z}\theta^*_+|^2d\textrm{vol}_g(x)dz\geq K\right\}\right|\leq\epsilon h/4.\]
Along the proof several smallness assumptions will be imposed on $\epsilon$, being a finite number this causes no problem for the argument to work. Notice that if $t\in I(h)\cap\mathcal{T}$ is such that $|\mathcal{A}(t)|\geq\frac{1}{4}|B^*_h|$ then using De Giorgi's as above we obtain
\[\int_{I(h)}\int_{B^*_g(h)}|\nabla_{x,z}\theta^*_+|^2d\textrm{vol}_g(x)dzdt\leq Ch^{n}\]
which follows from the local energy inequality, cf. \S\ref{localenergy}. Therefore one gets
\[|\mathcal{B}(t)|\leq\epsilon^{\frac{5}{2}}C h^{n+1}\]
using De Giorgi's isoperimetric inequality. This leads to $\|\theta\|_{L^2(Q_g(h))}\leq\epsilon h^{\frac{n+1}{2}}$ and $\|\theta^*\|_{L^2(Q^*_g(h))}\leq\epsilon h^{\frac{n}{2}+1}$. We claim that such estimate is true for some $t_1\in\mathcal{T}$ which does not lie on 
\[\left[-\frac{h}{4}+t^*,t^*\right]\subseteq [t^*-h,t^*]=I(h).\] 
Indeed, the non existence of $t^*$ leads to a contradiction with the size condition in the statement. We will use this facts to prove the following

\begin{claim}
The inequality $|\mathcal{A}(t)|\geq\frac{\epsilon^2}{2}|B_g^*(h)|$ holds for any $t\in I(h/4)\cap \mathcal{T}$. 
\end{claim}

The local energy inequality assures that for $t\geq t_1$ one has
\[\int_{B_g(h)}\theta_+^2(t)d\textrm{vol}_g(x)\leq\int_{B_g(h)}\theta_+^2(t_1)d\textrm{vol}_g(x)+C\|\theta\|^2_{L^{\infty}}(t-t_1)(h^{n-1}+h^n).\]
Observe that $t-t_1$ is of order $h$; let $t-t_1\leq\epsilon h$, if $\epsilon$ is small enough
\[\|\theta_+\|^2_{L^2(B_g(h))}\leq\frac{1}{100}h^n.\]
 Now, for such $t$
\[\begin{split}
\theta^*_+(x,t,z)&=\theta_+(x,t)+\int_0^z\partial_z\theta_+^*d\bar{z}\\
&\leq|\theta_+(x,0)|+\sqrt{z}\left(\int_0^z|\partial_z\theta_+^*|^2d\bar{z}\right)^{\frac{1}{2}}.
\end{split}\]
Applying the above with $z\leq\epsilon^2h$ implies
\[\|\theta^*_+(x,t,z)\|_{L^2(B_g(h))}\leq\|\theta_+(x,t)\|_{L^2(B_g(h))}+\epsilon h^{\frac{1}{2}}\left(\int_{B_g^*(h)}|\nabla_z\theta^*(x,t,\bar{z})|^2d\bar{z}d\textrm{vol}_g(x)\right)^{\frac{1}{2}}.\]
Our previous results show that it is bounded by $\frac{1}{100}h^n+\epsilon Ch^{\frac{n}{2}}$ which is smaller than $\frac{1}{2}|B_g(h)|^{\frac{1}{2}}$ if $\epsilon$ is small enough. Application of a weak $L^2$ inequality implies
\[\left|\left\{x\in B_g(h):\theta_+^*(x,t,z)\geq\frac{1}{2}\|\theta^*\|_{L^{\infty}(B_g^*(h))}\right\}\right|\leq\frac{1}{4}|B_g(h)|\]
Integrating $z\in I(\epsilon^2h)$ yields
\[\left|\left\{x\in B_g(h), z\in I(\epsilon^2h):\theta_+^*(x,t,z)\geq\frac{1}{2}\|\theta^*\|_{L^{\infty}(B_g^*(h))}\right\}\right|\leq\frac{\epsilon^2}{4}|B^*_g(h)|\]
Since $t\in\mathcal{T}$ one has that $|\mathcal{C}(t)|\leq 2\epsilon^6h^{n+1}\leq\epsilon^5|B_g^*(h)|$ and we get the following estimate from below
\[\begin{split}
|\mathcal{A}(t)|&\geq|B_g(h)\times I(\epsilon^2h)|-|\mathcal{C}(t)|\\
&\qquad-\left|\left\{x\in B_g(h), z\in I(\epsilon^2h):\theta_+^*(x,t,z)\geq\frac{1}{2}\|\theta^*\|_{L^{\infty}(B_g^*(h))}\right\}\right|\\
&\geq\left(\epsilon^2\left(1-\frac{1}{4}\right)-\epsilon^5\right)|B_g^*(h)|\geq\frac{1}{2}\epsilon^2|B_g^*(h)|
\end{split}\]
proving our claim, provided that $\epsilon$ is small enough.

Repeating the argument at the beginning of the proof one gets $|\mathcal{B}(t)|\leq C\sqrt{\epsilon}h^{n+1}$. So, finally, we have proved that, in fact, $|\mathcal{A}(t)|\geq\frac{1}{4}|B_g^*(h)|$ for any $t\in\mathcal{T}$ satisfying $t-t_1\leq\epsilon h$. Then, as before:
\[\begin{split}
|\mathcal{A}(t)|&\geq|B^*_g(h)|-|\mathcal{B}(t)|-|\mathcal{C}(t)|\\
&\geq\left(1-C\sqrt{\epsilon}-\epsilon^5\right)|B_g^*(h)|\geq\frac{1}{4}|B^*_g(h)|
\end{split}\]
which holds provided $\epsilon$ is small enough. This shows that $|\mathcal{A}(t)|\geq\frac{1}{4}|B_g^*(h)|$ holds in $t_1+I(\epsilon h)\cap\mathcal{T}$. We may change $t_1$ to some $t_2$ to its left provided $\mathcal{T}^c$ is small compared to $\epsilon h$, as it actually happens. Proceeding in this way one covers $\mathcal{T}\cap I(h/4)$. As a consequence
\[|\mathcal{B}(t)|\leq\epsilon^2|B^*_g(h)|\]
holds for any $t\in\mathcal{T}\cap I(h/4)$. And this allows us to estimate
\[\begin{split}
\int_{Q^*_g(h)}\left(\theta^*-\frac{1}{2}\|\theta^*\|_{L^{\infty}}\right)^2_+d\textrm{vol}_g(x)dzdt&=\|\theta^*\|^2_{L^{\infty}}\left(\frac{1}{2}\int_{t\in\mathcal{T}}|\mathcal{B}(t)|dt+\int_{t\notin\mathcal{T}}1dt\right)\\
&\leq C(\epsilon^2+\epsilon)|Q_g^*(h)|.
\end{split}\]
To get the other smallness condition one uses equation (\ref{isoperometrica2}), with $(\theta^*-1)_+$ instead, so that
\[\begin{split}
\int_{Q_g(h)}&\left(\theta-\frac{1}{2}\|\theta^*\|_{L^{\infty}}\right)_+^2dxdt=\int_{Q_g(h)}\left(\theta^*-\frac{1}{2}\|\theta^*\|_{L^{\infty}}\right)_+^2d\textrm{vol}_g(x)dt\\
&\quad-2\int_0^{h}\int_{0}^z\int_{B_g(h)}\left(\theta^*-\frac{1}{2}\|\theta^*\|_{L^{\infty}}\right)_+\partial_z\left(\theta^*-\frac{1}{2}\|\theta^*\|_{L^{\infty}}\right)_+d\textrm{vol}_g(x)d\bar{z}dt
\end{split}\]
But from the previous argument we know already that there exist some $z\in I(h)$ such that the first integral is $\sqrt{\epsilon}h^{n+1}$. By Cauchy-Schwarz inequality the other term is dominated by \[\|(\theta^*-\frac{1}{2}\|\theta^*\|_{L^{\infty}})_+\|_{L^2(Q_g^*(h))}\|\partial_z(\theta^*-\frac{1}{2}\|\theta^*\|_{L^{\infty}})_+\|_{L^2(Q_g^*(h))}.\] Finally, the last term can be bounded by the $L^2$ norm of the gradient $\nabla_{x,z}\theta^*$ which is controlled by $h^{\frac{n-1}{2}}$. All this altogether proves 
\[\left\|\left(\theta-\frac{1}{2}\|\theta^*\|_{L^{\infty}}\right)\right\|_{L^2(Q_g(h))}\leq\sqrt{\epsilon}h^{\frac{n+1}{2}}\]
which making $\epsilon$ even smaller, if necesary, implies Lemma ~\ref{isononlocal}.

Let us now explore the consequences of this rather technical Lemma:

\begin{proposition}\label{propbig}
For $h$ small enough, there exist a $\gamma<1$ (independent of the scale $h$) such that the following holds
\[\|\theta_+\|_{L^{\infty}(Q_g(h))}\leq\gamma\|\theta^*\|_{L^{\infty}(Q^*_g(2h))}.\]
\end{proposition}

\textsc{Proof:} without loss of generality we will assume that $\theta^*$ is such that $C_0=-\inf_{Q^*_g(h)}\theta^*=\sup_{Q^*_g(h)}\theta^*$ and such that $\theta$ is negative at least for half the points in $Q^*_g(h)$, otherwise one may substract an approppiate quantity or argue with $-\theta$ instead. We define the following truncations
\[\tau_k=2^k\left(\theta-\frac{1}{2}(1-2^{-k})C_0\right)_+.\]
Notice that the extension $\tau_k^*$ is precisely $2^k(\theta^*-\frac{1}{2}(1-2^{-k})C_0)_+$ and that all of them are, by construction, bounded by the same $C_0$. We claim now that $\tau_{k_0}$ satisfies the hypothesis of lemma \ref{isononlocal} for some integer $k_0\leq\frac{1}{\delta(\epsilon)}$. Indeed, otherwise
\[\{\tau^*_k<0\}=\left\{\tau^*_{k-1}<\frac{1}{2}C_0\right\}\geq\{\tau_{k-1}^*<0\}+\delta|Q_g(h)|\]
which can not hold inductively longer than $\frac{1}{\delta}$ times. As a consequence, Lemma \ref{isononlocal} implies that $\tau_{k_0}$ is under the hypothesis of Proposition \ref{propsmall}. Hence
\[\|(\tau_{k_0})_+\|_{L^{\infty}(Q_g(h/2))}\leq\gamma\|\tau_{k_0}^*\|_{L^{\infty}(Q^*_g(h))}\]
for some $\gamma<1$. Unravelling notation one observes that
\[\|\theta_+\|_{L^{\infty}(Q_g(h))}\leq \left(1-\frac{1-\gamma}{2^{k_0+1}}\right)\|\theta^*\|_{L^{\infty}(Q^*_g(2h))}.\]
Notice now that the decrease is smaller than one and that it can be taken independently of the scale $h$ we are working with.

\subsection{Proof of Theorem \ref{II}}

\begin{figure}[htb]
\centering
\includegraphics[width=90mm]{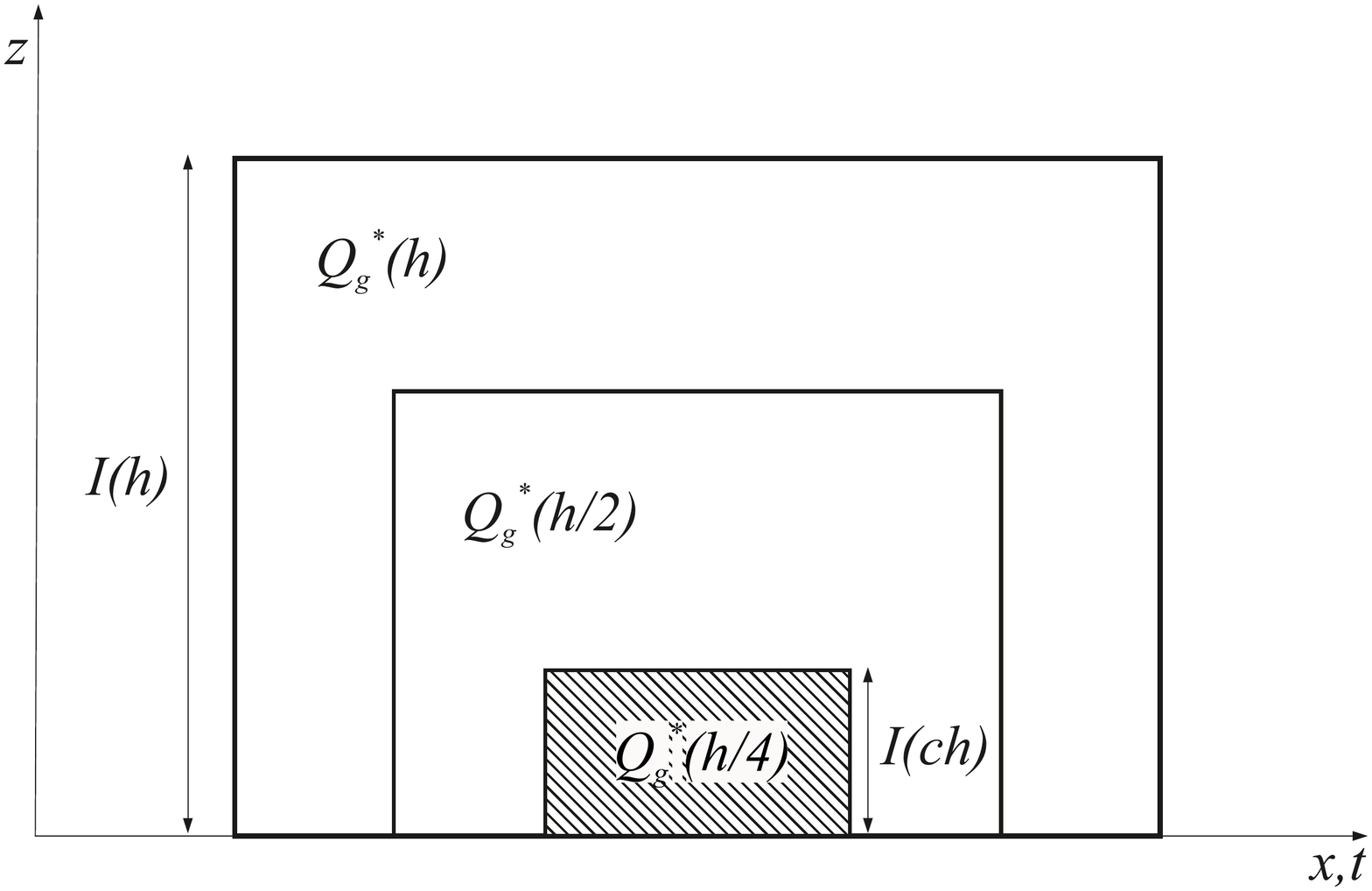}
\caption{}\label{fig:b4}
\end{figure}

Using again a barrier of the form 
\[\|\theta^*\|_{L^{\infty}(Q^*_g(h))}b_2+L(z)\]
where $L(z)$ is linear function interpolating between $L(0)=\|\theta\|_{L^{\infty}(Q_g(h/2))}$ and $L(h/2)=\|\theta^*\|_{L^{\infty}(Q^*_g(h)}$ (see figure \ref{fig:b4}). As a consequence (cf. Lemma \ref{b2}) one may show the existence of $\gamma<1$ such that
\[\|\theta^*\|_{L^{\infty}(Q^*_g(ch))}\leq\gamma\|\theta^*\|_{L^{\infty}(Q^*_g(h))}\]
for some sufficiently small positive constant $c$. This implies the statement with $\alpha=\alpha(c,\|\theta_0\|,M,t_0)>0$ by an argument analogous to the one provided in \S\ref{cafvas}.

\section{Proof of theorem \ref{I}} \label{CVproof}

In this section we provide a variant of an argument from \cite{CV} where the authors use translations and dilations to construct a sequence of $F_k$ related to $\theta$ for which the oscillation decays fast enought to obtain $C^{\alpha}$ regularity.  

Let us begin illustrating a geometrical idea that we take advantage of and providing a glimpse of the method: first we define an auxiliary function
\[F(s,x)=\theta(R_s(x),s+t_0)\]
where $R_s$ is a rigid rotation of the sphere around some axis and $s$ denotes the arc length of the particles moving in the corresponding equator. This device compensates the probably high velocities of $u$. Notice that such an $F$ satisfies the same equation
\[\partial_sF+v\cdot\nabla_gF=-\Lambda F\]
where $R$ is the infinitesimal generator of the rotation $R_s$ and $v=u-\dot{R}_s$. We used rotations because they are global isometries respecting the nonlocal diffusive operator. 

First we need to check that the velocity satisfies the hypothesis 
\begin{equation}
\int_{B_g(h)}|v(x,t)|^{2n}d\textrm{vol}_g(x)\leq C(M)h^n\label{drifthip}\end{equation}
uniformly on time. This will be done locally of course. Fix some $x_0\in\mathbb{S}^2$ and a ball $B_g(h)$ around it. 


We use the standard embedding $\mathbb{S}^2\subseteq\mathbb{R}^3$. Let $u(x)=u_t(x)+u_n(x)$ where $u_t$ denotes the projection to the tangent plane at $x_0$, $T_{x_0}\mathbb{S}^2$. Now, near $x_0$ it is evident that $u_n$ is small, more precisely
\[u_n(x)=O(ru)\textrm{ for $x\in B_g(r)$}\]
provided $r$ is small enough. We define $R_s$ as the rotation generated by the following tangent vector at $x_0$
\[\dot{R}_s(x_0)=\frac{1}{|B_g(h)|}\int_{B_g(h)}u_t(s+t,R_s(x))d\textrm{vol}_g(x)\in T_{x_0}\mathbb{S}^2.\]
This definition is equivalent to an ordinary differential equation for $R_s$ with $R_0=\textrm{id}$. Notice that, for the same reason as above,
\[R(x)=\dot{R}_s(x_0)+O(r\dot{R}_s(x_0))\textrm{ for any $x\in B_g(r)$}\]
for small $r$. Hence
\[\dot{R}_s(x_0)=\frac{1}{|B_g(h)|}\int_{B_g(h)}u \ d\textrm{vol}_g(x)-\frac{1}{|B_g(h)|}\int_{B_g(h)}u_n \ d\textrm{vol}_g(x)=u_{B_g(h)}+O(1)\]
where
\[u_{B_{g}(h)}:=\frac{1}{|B_g(h)|}\int_{B_g(h)}u \ d\textrm{vol}_g(x).\]
Next we use Cauchy-Schwarz to bound the second summand, $L^2$ boundedness of the Riesz transform and the extra $h$ from the estimate above:
\[\begin{split}
\frac{1}{|B_g(h)|}\int_{B_g(h)}u_n \ d\textrm{vol}_g(x)&\leq \frac{1}{|B_g(h)|^{1/2}}\left(\int_{B_g(h)}|u_n|^2 \ d\textrm{vol}_g(x)\right)^{1/2}\\
&\lesssim h^{1-n/2}\left(\int_{B_g(h)}|u|^2 \ d\textrm{vol}_g(x)\right)^{1/2}=O(1)
\end{split}\]
This is where the two dimensionality of the sphere becomes crucial in the argument. For higher dimensional spheres the bound is not good enough. We can now estimate the $L^{2n}$ norm
\[\begin{split}
\left(\frac{1}{|B_g(h)|}\int_{B_g(h)}|u-\dot{R}_s|^{2n}d\textrm{vol}_g(x)\right)^{1/2n}&\leq\left(\frac{1}{|B_g(h)|}\int_{B_g(h)}|u-u_{B_g(h)}|^{2n}d\textrm{vol}_g(x)\right)^{1/2n}\\
&\quad+\left(\frac{1}{|B_g(h)|}\int_{B_g(h)}|u_{B_g(h)}-\dot{R}_s(x_0)|^{2n}d\textrm{vol}_g(x)\right)^{1/2n}\\
&\qquad+\left(\frac{1}{|B_g(h)|}\int_{B_g(h)}|\dot{R}_s(x_0)-\dot{R}_s|^{2n}d\textrm{vol}_g(x)\right)^{1/2n}
\end{split}\]
the first is bounded due to the John-Nirenberg inequality since $u\in BMO(\mathbb{S}^2,\mathbb{R}^3)$ (cf. \cite{BN, Tay}). The second is $O(1)$ from the above estimates while the third is bounded by $h|\dot{R}_s(x_0)|$, which is again $O(1)$ by Cauchy-Schwarz inequality due to the extra $h$. This shows that hypothesis ~\ref{drifthip}
\[\int_{B_g(h)}|v|^{2n}d\textrm{vol}_g(x)\leq Ch^n\]
is satisfied for $n=2$. 
Now it is licit to infer from Theorem \ref{II} the existence of some $\gamma<1$ for which
\[\textrm{osc}_{Q_g(h)}F\leq\gamma\cdot\textrm{osc}_{Q_g(2h)}F\]
Notice that this can not be rephrased directly in terms of the oscillation of $\theta$ due to the possible high velocities that may displace points too far away. One might, nevertheless, bound it from below paying the price of making times smaller to compensate this high velocities while one may use the a priori bound for the displacement
\[R_s(x)\leq s\sup_{s\leq h} \dot{R}_s(x)\leq Csh^{-n}\sup_{s\leq h}\int_{B_g(h)}u\leq C h^{n/p'-n+1}\|u\|_{L^p(M)}\]
being $u$ a Riesz transform which is uniformly bounded. Therefore, choosing $p>n$ we can estimate it from above by, say, $h^{1/3}$. As a consequence we get an inequality of the form
\[\textrm{osc}_{Q_g(h^{K})}\theta\leq\gamma\cdot\textrm{osc}_{Q_g(h)}\theta\]
for some $K$ big enough and $h$ small enough (depending on fixed quantities). This implies a  modulus of continuity of the form $\omega(\rho)=\log(1/\rho)^{-\alpha}$ for some $\alpha=\alpha(t_0,\|\theta_0\|_{L^2(M)})$ which deteriorates as $t_0$ approaches zero.

\section{\textbf{Appendix: global weak solutions}}\label{debiles}

In this section we provide the existence theory of global weak solutions to incompressible drift diffusion equations on compact manifolds given by
\[\left\{\begin{array}{l}
\partial_t\theta+u\cdot\nabla_g\theta+\Lambda\theta=0 \\
\nabla\cdot u= 0
\end{array}\right.\]
where the velocity field $u=\Psi[\theta]$, with $\Psi$ a zero order pseudodifferential operator. More precisely, we prove the next result:
\begin{theorem}\label{weak} Let $\theta_{0}\in L^{2}(M)$ be the initial data and let $T>0$. Then there exists a weak solution $\theta \in L^{\infty}(0,T;L^{2}(M))\cap L^{2}(0,T;{H}^{1/2}(M))$ to the equation above i.e, the following equality holds
\[ \int_{M} \theta \psi \ d\textrm{vol}_{g}(x) +\int_{0}^{T} \int_{M} \theta\Lambda\psi-\theta u\cdot \nabla_{g}\psi d\textrm{vol}_{g}(x) \ dt= \int_{M} \theta_{0} \psi \ d\textrm{vol}_{g}(x)  \]
for each test function $\psi\in C^{\infty}([0,T]\times M)$.
\end{theorem}

\textsc{Proof of Theorem ~\ref{weak}:} we will prove it using Galerkin approximations. For each integer $m>0$, consider the truncations $\theta_{m}$ given by
\[ \theta_{m}(x,t)= \displaystyle\sum_{k=0}^{m} f_{k}(t)Y_{k}(x),\]
where $Y_{k}$ are the eigenfunctions of $-\Delta_{g}$ with eigenvalues $\lambda^{2}_{k}$. Denote by $\mathbb{P}_{m}$ the projection onto the space generated by $\lbrace Y_{k}\rbrace_{k=1}^{m}$. Then for each fixed $m$, let us consider the truncated system:
\[\left\{\begin{array}{l}
\partial_{t} \theta_{m}+\mathbb{P}_{m}(u_{m}\cdot \nabla_{g}\theta_{m}) + \Lambda \theta_{m}=0 \\
u_{m}= \Psi[\theta_{m}] \\
\theta_{m}(x,0)= \mathbb{P}_{m}\theta_{0}(x) 
\end{array}\right.\]
Although the above system seems to be a partial differential equation, it can be interpreted as a system of ordinary differential equations for the coefficients, $f_{k}(t)$,
\[ f'_{k}(t)=-\lambda_{k}f_{k}(t) + \displaystyle\sum_{\ell=1}^{m} a_{k,l}(t)f_{\ell}(t),\]
where $a_{k,\ell}(t)=\int_{M} u(x,t) \cdot \nabla_{g}Y_{k}(x)Y_{\ell}(x) \ d\textrm{vol}_{g}(x)$ and the initial condition is provided by $\theta_m(0)=\mathbb{P}_m(\theta_0)$. The initial condition has bounded energy which, taking into account the nature of $\Psi$, implies the velocity is bounded in $L^2$. From this it is easily proved that the coefficients $a_{k,\ell}$ are uniformly bounded. This allows us to use standard Picard-Lindel{\"o}f existence and uniqueness theorem for ordinary differential equations to show global existence of $f_{k}(t)$ (cf. \cite{CL}, p. 20). Moreover, thanks to the divergence free condition the nonlinear term vanishes after integration by parts
\[ \int_{M} \mathbb{P}_{m}(u_{m}\cdot\nabla_{g}\theta_{m}) \theta_{m} d\textrm{vol}_{g}(x)=0\]
and we get the uniform energy estimate
\[ \|\theta_{m}(T)\|^{2}_{L^{2}(M)} + \int_{0}^{T}\|\theta_{m}(\tau)\|^{2}_{{H}^{1/2}(M)} d\tau\leq \|\theta_{0}\|^{2}_{L^{2}(M)}.\] 
This implies $\theta_m$ is uniformly bounded in $L^{\infty}(0,T;L^{2}(M))\cap L^{2}(0,T;H^{1/2}(M))$ for every $T>0$. This ensures that, up to subsequence, $\theta_{m}$ converges weakly to some $\theta\in  L^{2}(0,T;H^{1/2}(M))$ the limit still belongs to $L^{\infty}(0,T;L^2(M))$. We may force the subsequence to be such that $\theta_m(T)$ converges weakly to $\theta(T)$ in $L^2(M)$. In what remains we will show that the limit function $\theta$ obtained above is a weak solution of the Cauchy initial problem as stated. 

Testing the truncated equation against some $\psi\in C^{\infty}([0,T]\times M)$  we obtain
\[\begin{split}
\int_{M}\theta_{m}(x,T)\psi \ d\textrm{vol}_{g}(x) &+ \int_{0}^{T} \int_{M} \theta_{m}\Lambda\psi d\textrm{vol}_{g}(x)\\
& - \int_{0}^{T}\int_{M} (\theta_{m}u_{m})\cdot\nabla_{g}\psi d\textrm{vol}_{g}(x) dt = \int_{M}P_{m}\theta_{0}\psi d\textrm{vol}_{g}(x).
\end{split}\]
If we can take limits inside the integrals the weak formulation would be satisfied and we would be done. Since $\theta_m(T)$ converges weakly to $\theta(T)$ in $L^2$ it follows that
\[\int_{M} \theta_m(x,T)\psi \ d\textrm{vol}_{g}(x)\rightarrow \int_{M} \theta(x,T)\psi \ d\textrm{vol}_{g}(x)\]
while since $\|\mathbb{P}_{m}f-f\|_{L^{2}(M)}\rightarrow 0$ as we tend $m\rightarrow \infty$ for any $f\in L^2(M)$ it follows that
\[\int_{M} \mathbb{P}_{m}\theta_{0}(x)\psi \ d\textrm{vol}_{g}(x)\rightarrow \int_{M} \theta_{0}(x)\psi \ d\textrm{vol}_{g}(x).\]
The convergence of the second term above follows from the weak convergence in $L^2(0,T;H^{1/2}(M))$, indeed
\[ \int_{0}^{T}\int_{M} \theta_{m}\Lambda\psi  d\textrm{vol}_{g}(x) \rightarrow \int_{0}^{T}\int_{M} \theta\Lambda\psi  d\textrm{vol}_{g}(x) \]
since $\Lambda\psi$ is still a test function. Regarding the nonlinear term, we realize that the weak convergence is not sufficient. We need some stronger convergence. It would be enough to prove that
\[ \theta_{m}\rightarrow \theta \ \text{strongly in} \ L^{2}(0,T;L^{2}(M))\] 
Let us show how this actually helps to deal with the nonlinear term. Indeed, for any test $\psi\in C^{\infty}([0,T]\times M)$ we would have
\[\begin{split}
\left|\int_{0}^{T}\int_{M}\theta_{m} u_m\cdot \nabla_{g}\psi-\theta u\cdot\nabla_{g}\psi  d\textrm{vol}_{g}(x)dt \right| &\leq \left|\int_{0}^{T}\int_{M}\left(\theta_{m}-\theta\right)u_{m}\cdot \nabla_{g}\psi  d\textrm{vol}_{g}(x)dt \right|\\
&\quad\quad+ \left|\int_{0}^{T}\int_{M}\theta\left(u_{m}-u\right)\cdot \nabla_{g}\psi d\textrm{vol}_{g}(x)dt \right|
\end{split}\]
which can be bounded by 
\[\begin{split}
C(\psi)&\left\{\|\theta_{m}-\theta\|_{L^2(0,T;L^{2}(M))}\|u_{m}\|_{L^2(0,T;L^{2}(M))}\right.\\
&\qquad\qquad\qquad\quad\left.+\|u_{m}-u\|_{L^2(0,T;L^{2}(M))}\|\theta\|_{L^2(0,T;L^{2}(M))}\right\},
\end{split}\]
and now it becomes obvious that it converges to zero  as $m\rightarrow\infty$.

To prove the claimed strong convergence we will proced through the Aubin-Lions compactness lemma. To do so let us first observe that
\[ \partial_{t}\theta_{m} \  \text{is bounded in} \ L^{(n+1)/n}(0,T;B) \] 
where $B=W^{-1,\frac{n+1}{n}}(M)+H^{-1/2}(M)$. Indeed the energy estimate above implies (see \S3 for details) that 
\[\theta_{m}\textrm{ is uniformly bounded in }L^{2(n+1)/n}(0,T;L^{2(n+1)/n}(M)).\]
 As before, recall $\Psi$ is a zero order operator it conserves $L^{p}(M)$ norms for $p\in(1,\infty)$, it follows then that
 \[u_{m}\textrm{ is also bounded in }L^{2(n+1)/n}(0,T;L^{2(n+1)/n}(M)).\]
 Therefore, $\nabla_g\cdot(u_{m}\theta_{m})$ is bounded in $L^{(n+1)/n}(0,T; W^{-1,(n+1)/n}(M))$. Lastly, $\Lambda\theta_{m}$ is bounded in $L^{2}(0,T; H^{-1/2}(M))$. From the equation itself we get that $\partial_{t}\theta_{m}$ is bounded in $L^{(n+1)/n}(0,T; B)$. Now $H^{1/2}(M)\rightarrow L^2(M)$ is a compact inclusion by Rellich's theorem while by Sobolev $L^2(M)$ is continuously embbedded into $B$. As a consequence the Aubin-Lions lemma provides the strong convergence (cf. \cite{JL}). Finally, notice that $T>0$ was arbitrary, proving that $\theta$ is a global weak solution.

\section{\textbf{Acknowledgments}}

The authors are indebted to T. Pernas-Casta\~no who carefully draw the figures herein.

This work has been partially supported by ICMAT-Severo Ochoa project SEV-2015-0554 and the MTM2011-2281 project of the MCINN (Spain).


\begin{thebibliography}{10}
\bibitem{AOMC}
Alonso-Or\'an, D.; C\'ordoba, A.; Mart\'inez, A. D., {\em Integral representation for the Laplace-Beltrami operators}, unpublished.
\bibitem{Au}
Aubin, T., {\em Nonlinear Analysis on Manifolds. Monge-Amp\`ere equations}, Springer 1982.
\bibitem{BC}
Balodis, P.; C\'ordoba, A., {\em An inequality for Riesz transforms implying blow-up for some nonlinear and nonlocal transport equations}, Adv. of Math. Vol. 214 (10) (2007), pp. 1-39.
\bibitem{BN}
Brezis, H.; Nirenberg, L., {\em Degree theory and BMO, part I: compact manifolds without boundaries}, Selecta Mathematica 1 (1995), pp. 197-263.
\bibitem{CCW}
Constantin, P.; C\'ordoba, D.; Wu, J., {\em On the critical dissipative quasi-geostrophic equation}, Indiana Unif. Math. J. 50 (2001), pp. 97-107.
\bibitem{CSil} Caffarelli, L. A.; Silvestre, L., {\em An extension problem related to the fractional Laplacian}, Communications in Partial Differential Equations, 32 (2007) Issue 8, 1245.
\bibitem{CS}
Caffarelli, L. A.; Sire, Y., {\em On some pointwise inequalities involving nonlocal operators}, Arxiv 1604.05665.
\bibitem{CV} 
Caffarelli, L. A.; Vasseur, A., {\em Drift diffusion equations with fractional diffusion and the quasi-geostrophic equation}, Ann. of Math. Vol. 171, No. 3 (2010), pp. 1903-1930.
\bibitem{CV2}
Caffarelli, L. A.; Vasseur, A., {\em The De Giorgi Method for Nonlocal Fluid Dynamics}.
\bibitem{Ch}
Chavel, I., {\em Eigenvalues in Riemannian geometry}, Academic Press, 1984.
\bibitem{Ch2}
Chavel, I, {\em Riemannian geometry: a modern introduction}, Cambridge University Press, 1993.
\bibitem{CL}
Coddington, E. A.; Levinson, N., {\em Theory of Ordinary and Differential Equations}, McGraw-Hill, 1955.
\bibitem{CI1}
Constantin, P.; Ignatova, M., {\em Critical SQG in bounded domains}, to appear in Annals of PDE, 2016.
\bibitem{CI2}
Constantin, P.; Ignatova, M., {\em Nonlinear lower bounds for the fractional Laplacian with Dirichlet boundary conditions and applications},  to appear in Int. Math. Res. Notices, 2016.
\bibitem{CMT}
Constantin, P.; Majda, A.; Tabak, E., {\em Formation of strong fronts in the 2D quasi-geostrophic thermal active scalar}, Nonlinearity 7 (1994), pp. 1495-1533.
\bibitem{CVi}
Constantin, P.; Vicol, V., {\em Nonlinear maximum principles for dissipative linear nonlocal operators and applications}, Geometric and Functional Analysis, 22 (2012), no. 5, 1289-1321. 
\bibitem{CW}
Constantin, P.; Wu, J., {\em Behaviour of solutions of 2D quasigestrophic equations}, SIAM J. Math. Anal. 30(5) (1999), pp. 937-948. 
\bibitem{Cdie}
C\'ordoba, D., {\em Nonexistence of simple hyperbolic blow-up for the quasi-geostrophic equation}, Ann. of Math. 148 (1998), pp. 1135-1152.
\bibitem{CC} 
C\'ordoba, A.; C\'ordoba, D., {\em A Maximum Principle Applied to Quasi-Geostrophic Equations}, Commun. Math. Phys. 249 (2004), pp. 511-528.
\bibitem{CC2} 
C\'ordoba, A.; C\'ordoba, D., {\em A pointwise estimate for fractionary derivatives with applications to partial differential equations.} Proceedings of the National Academy of Sciences of the United States of America 100 (26) (2003), pp. 15316-15317.
\bibitem{CCF}
C\'ordoba, A.; C\'ordoba, D.; Fontelos, M. A., {\em Formation of singularities for a transport equation with nonlocal velocity}, Ann. of Math., 162 (2005), pp. 1377-1389.
\bibitem{CCG} 
C\'ordoba, A.; C\'ordoba, D.; Gancedo, F., {\em Interface evolution: the Hele-Shaw and Muskat problems}, Ann. of Math. Vol. 173 (2011), pp. 477-542.
\bibitem{CM}
C\'ordoba, A.; Mart\'inez, A. D., {\em A pointwise inequality for fractional laplacians}, Adv. of Math., Vol. 280 (2015), pp. 79-85.
\bibitem{DG}
De Giorgi, E., {\em Sulla differenziabilit\`a e l'analiticit\`a delle estremali degli integrali multipli regolari}, Mem. Accad. Sci. Torino Cl. Sci. Fis. Mat. Nat. (3) 3 (1957), pp. 25-43.
\bibitem{DZ}
Duandikoetxea, J., {\em An\'alisis de Fourier}, Ediciones UAM, 1989.
\bibitem{held}
Held, I. M.; Pierrehumbert, R. T.; Garner, S. T.; Swanson, K. L., {\em Surface quasi-geostrophic dynamics}, J. Fluid Mech. 282 (1995), pp. 1-20.
\bibitem{LE}
Evans, L. C., {\em Partial differential equations}, Graduate Studies in Mathematics Vol. 19, 1998.
\bibitem{Fo}
Folland, G. B., {\em Introduction to Partial Differential Equations}, second edition. Princeton University Press, 1995.
\bibitem{Ho}
H\"ormander, L., {\em The analysis of linear partial differential operators (vol. III)}, Grundlehren der mathematischen Wissenschaften, Vol. 274, Springer, 1985.
\bibitem{Ke}
Kellogg, O. D., {\em Foundations of Potential Theory}, Dover, 1929.
\bibitem{KNV}
Kiselev, A.; Nazarov, F.; Volberg, A., {\em Global well-posedness for the critical 2D dissipative quasigeostrophic equation}, Invent. Math. 167 (2007), pp. 445-453.
\bibitem{JL}
Lions, J.L., {\em Quelque M\'ethodes de R\'esolutions des Probl\`emes aux Limites Non-Lin\'eaires}, Dunod, Paris, 1969.
\bibitem{N}
Nash, J., {\em Continuity of Solutions of Parabolic and Elliptic Equations}, American Journal of Mathematics
Vol. 80, No. 4 (1958), pp. 931-954.
\bibitem{PW}
Protter, M. H.; Weinberger, H. F., {\em Maximum principles in differential equations}, Pretice-Hall, 1967.
\bibitem{R}
Resnick, S. G. {\em Dynamical problems in non-linear advective partial differential equations}, PhD thesis (Chicago University, 1995).
\bibitem{So}
Sogge, Ch. D., {\em The Hangzhou lectures on Eigenfunctions of the Laplacian}, Princeton University Press, 2014.
\bibitem{S0}
Stein, E. M., {\em Harmonic Analysis: Real-Variable Methods, Orthogonality and Oscillatory Integrals}, Princeton University Press, 1993.
\bibitem{S1}
Stein, E. M., {\em Singular integrals and Differentiability Properties of Functions}, Princeton University Press, 1970.
\bibitem{Tay}
Taylor, M., {\em Hardy Spaces and bmo on Manifolds with Bounded Geometry}, J. Geom. Anal. 19 (1) (2009), pp. 137-190.
\bibitem{V}
Vasseur, A. F., {\em The De Giorgi method for elliptic and parabolic equations and some applications}, to appear in Lectures on the Analysis of Nonlinear Partial Differential Equations Vol. 4. 
\bibitem{W}
Widman, K., {\em Inequalities for the Green function and boundary continuity of the gradient of solutions of elliptic differential equations}, Math. Scan. 21 (1967), pp. 17-37.
\end{thebibliography}
\end{document}